# Sample size and positive false discovery rate control for multiple testing

Zhiyi Chi*

*Department of Statistics*
*University of Connecticut*
*215 Glenbrook Road, U-4120*
*Storrs, CT 06269*
*e-mail:* zchi@stat.uconn.edu

**Abstract:** The positive false discovery rate (pFDR) is a useful overall measure of errors for multiple hypothesis testing, especially when the underlying goal is to attain one or more discoveries. Control of pFDR critically depends on how much evidence is available from data to distinguish between false and true nulls. Oftentimes, as many aspects of the data distributions are unknown, one may not be able to obtain strong enough evidence from the data for pFDR control. This raises the question as to how much data are needed to attain a target pFDR level. We study the asymptotics of the minimum number of observations per null for the pFDR control associated with multiple Studentized tests and $F$ tests, especially when the differences between false nulls and true nulls are small. For Studentized tests, we consider tests on shifts or other parameters associated with normal and general distributions. For $F$ tests, we also take into account the effect of the number of covariates in linear regression. The results show that in determining the minimum sample size per null for pFDR control, higher order statistical properties of data are important, and the number of covariates is important in tests to detect regression effects.

## 1. Introduction

A fundamental issue for multiple hypothesis testing is how to effectively control Type I errors, namely the errors of rejecting null hypotheses that are actually true. The False Discovery Rate (FDR) control has generated a lot of interest due to its more balanced trade-off between error rate control and power than the traditional Familywise Error Rate control (1). For recent progress on FDR control and its generalizations, see (6–12, 14–16, 19) and references therein.

Let $R$ be the number of rejected nulls and $V$ the number of rejected true nulls. By definition, $\text{FDR} = E[V/(R \vee 1)]$. Therefore, in FDR control, the case $R = 0$ is counted as "error-free", which turns out to be important for the controllability of the FDR. However, multiple testing procedures are often used in situations where one explicitly or implicitly aims to obtain a nonempty set of rejected nulls. To take into account this mind-set in multiple testing, it is appropriate to control the positive FDR (pFDR) as well, which is defined as $E[V/R \,|\, R > 0]$ (18). Clearly, when all the nulls are true, the pFDR is 1 and therefore cannot be controlled. This is a reason why the FDR is defined as it is (1). On the other hand, even when there is a positive proportion of nulls that are false, the pFDR can still be significantly greater than the FDR, such that when some nulls are indeed rejected, chance is that a large proportion or even almost all of them are falsely rejected (3, 4).

---

*Research partially supported by NIH grant MH68028.

The gap between FDR and pFDR arises when the test statistics cannot provide arbitrarily strong evidence against nulls (4). Such test statistics include $t$ and $F$ statistics (3). These two share a common feature, that is, they are used when the standard deviations of the normal distributions underlying the data are unknown. In reality, it is a rule rather than exception that data distributions are only known partially. This suggests that, when evaluating rejected nulls, it is necessary to realize that the FDR and pFDR can be quite different, especially when the former is low.

In order to increase the evidence against nulls, a guiding principle is to increase the number of observations for each null, denoted $n$ for the time being. In contrast to single hypothesis testing, for problems that involve a large number of nulls, even a small increase in $n$ will result in a significant increase in the demand on resources. For this reason, the issue of sample size per null for multiple testing needs to be dealt with more carefully. It is known that FDR and other types of error rates decrease in the order of $O(\sqrt{\log n/n})$ (13). In this work, we will consider the relationship between $n$ and pFDR control, in particular, for the case where false nulls are hard to separate from true ones. The basic question to be considered is: in order to attain a certain level of pFDR, what is the minimum value for $n$. This question involves several issues. First, how does the complexity of the null distribution affect $n$? Second, is normal or $t$ approximation appropriate in determining $n$? In other words, is it necessary to incorporate information on higher order moments of the data distribution? Third, what would be an attainable upper bound for the performance of a multiple testing procedure based on partial knowledge of the data distributions?

In the rest of the section, we first set up the framework for our discussion, and then outline the other sections.

### 1.1. Setup and basic approach

Most of the discussions will be made under a random effects model (10, 18). Each null $H_i$ is associated with a distribution $F_i$ and tested based on $\xi_i = \xi(X_{i1}, \ldots, X_{in})$, where $X_{i1}, \ldots, X_{in}$ are iid $\sim F_i$ and the function $\xi$ is the same for all $H_i$. Let $\theta_i = \mathbf{1}\,\{H_i \text{ is true}\}$. The random effects model assumes that

$(\theta_i, \xi_i)$ are independent, such that

$$\theta_i \sim \text{Bernoulli}(\pi), \quad \xi_i \mid \theta_i \sim \begin{cases} P_0^{(n)} \text{ with density } p_0^{(n)}, & \text{if } \theta_i = 0 \\ P_1^{(n)} \text{ with density } p_1^{(n)}, & \text{if } \theta_i = 1 \end{cases} \quad (1.1)$$

where $\pi \in [0, 1]$ is a fixed population proportion of false nulls among all the nulls. Note that $P_i^{(n)}$ of depend on $n$, the number of observations for each null. It follows that the minimum pFDR is (cf. (4))

$$\alpha_* = \frac{1-\pi}{1-\pi+\pi\rho_n}, \quad \text{with} \quad \rho_n := \sup \frac{p_1^{(n)}}{p_0^{(n)}}. \quad (1.2)$$

In order to attain pFDR $\leq \alpha$, there must be $\alpha_* \leq \alpha$, which is equivalent to $(1-\alpha)(1-\pi)/(\alpha\pi) \leq \rho_n$. For many tests, such as $t$ and $F$ tests, $\rho_n < \infty$ and $\rho_n \uparrow \infty$ as $n \to \infty$. Then, the minimum sample size per null is

$$n_* = \min\{n : (1-\alpha)(1-\pi)/(\alpha\pi) \leq \rho_n\}. \quad (1.3)$$

In general, the smaller the difference between the distributions $F_i$ under false nulls and those under true nulls, the smaller $\rho_n$ become, and hence the larger $n_*$ has to be. Our interest is how $n_*$ should grow as the difference between the distributions tends to 0.

**Notation** Because $(1-\alpha)(1-\pi)/(\alpha\pi)$ regularly appears in our results, it will be denoted by $Q_{\alpha,\pi}$ from now on.

### 1.2. Outlines of other sections

Section 2 considers $t$ tests for normal distributions. The nulls are $H_i : \mu_i = 0$ for $N(\mu_i, \sigma_i)$, with $\sigma_i$ unknown. It will be shown that if $\mu_i/\sigma_i \equiv r$ for false nulls, then, as $r \downarrow 0$, the minimum sample size per null $\sim (1/r) \ln Q_{\alpha,\pi}$ and therefore it depends on at least 3 factors: 1) the target pFDR control level, $\alpha$, 2) the proportion of false nulls among the nulls, $\pi$, 3) and the distributional properties of the data, as reflected by $\mu_i/\sigma_i$. In contrast, for FDR control, there is no constraint on the sample size per null. The case where $\mu_i/\sigma_i$ associated with false nulls are sampled from a distribution will be considered as well. This section also illustrates the basic technique used throughout the article.

Section 3 considers $F$ tests. The nulls are $H_i : \boldsymbol{\beta}_i = 0$ for $Y = \boldsymbol{\beta}_i^T \boldsymbol{X} + \epsilon$, where $\boldsymbol{X}$ consists of $p$ covariates and $\epsilon \sim N(0, \sigma_i)$ is independent of $\boldsymbol{X}$. Each $H_i$ is tested with the $F$ statistic of a sample $(Y_{ik}, \boldsymbol{X}_k)$, $k = 1, \ldots, n+p$, where $n \geq 1$ and $\boldsymbol{X}_1, \ldots, \boldsymbol{X}_{n+p}$ consist a fixed design for the nulls. Note that $n$ now stands for the difference between the sample size per null and the number of covariates included in the regression. The asymptotics of $n_*$, the minimum value for $n$ in order to attain a given pFDR level, will be considered as the regression effects become increasingly weak and/or as $p$ increases. It will be seen that $n_*$ must stay positive. The weaker the regression effects are, the larger $n_*$ has to be. Under certain conditions, $n_*$ should increase at least as fast as $p$.

Section 4 considers $t$ tests for arbitrary distributions. We consider the case where estimates of means and variances are derived from separate samples, which allows detailed analysis with currently available tools, in particular, uniform exact large deviations principle (LDP) (2). It will be shown that the minimum sample size per null depends on the cumulant generating functions of the distributions, and thus on their higher order moments. The asymptotic results will be illustrated with examples of uniform distributions and Gamma distributions. An example of normal distributions will also be given to show that the results are consistent with those in Section 2. We will also consider how to split the random samples for the estimation of mean and the estimation of variance in order to minimize the sample size per null.

Section 5 considers tests based on partial information on the data distributions. The study is part of an effort to address the following question: when knowledge about data distributions is incomplete and hence Studentized tests are used, what would be the attainable minimum sample size per null. Under the condition that the actual distributions belong to a parametric family which is unknown to the data analyzer, a Studentized likelihood test will be studied. We conjecture that the Studentized likelihood test attains the minimum sample size per null. Examples of normal distributions, Cauchy distributions, and Gamma distributions will be given.

Section 6 concludes the article with a brief summary. Most of the mathematical details are collected in the Appendix.

## 2. Multiple *t*-tests for normal distributions

### 2.1. Main results

Suppose we wish to conduct hypothesis tests for a large number of normal distributions $N(\mu_i, \sigma_i)$. However, neither $\sigma_i$ nor any possible relationships among $(\mu_i, \sigma_i)$, $i \geq 1$, are known. Under this circumstance, in order to test $H_i : \mu_i = 0$ simultaneously for all $N(\mu_i, \sigma_i)$, an appropriate approach is to use the $t$ statistics of iid samples $Y_{i1}, \ldots, Y_{i,n+1} \sim N(\mu_i, \sigma_i)$:

$$T_i = \frac{\sqrt{n+1}\,\bar{Y}_i}{S_i}, \quad \bar{Y}_i = \frac{1}{n+1}\sum_{j=1}^{n+1} Y_{ij}, \quad S_i^2 = \frac{1}{n}\sum_{j=1}^{n+1}(Y_{ij} - \bar{Y}_i)^2. \qquad (2.1)$$

Suppose the sample size $n+1$ is the same for all $H_i$ and the samples from different normal distributions are independent of each other.

Under the random effects model (1.1), we first consider a case where distributions with $\mu_i \neq 0$ share a common characteristic, i.e., signal-noise ratio defined in the remark following Theorem 2.1.

**Theorem 2.1.** *Under the above condition, suppose that, unknown to the data analyzer, when $H_i$ is false, $\mu_i/\sigma_i = r > 0$, where $r$ is a constant independent* of $i$. Given $0 < \alpha < 1$, let $n_*$ be the minimum value of $n$ in order to attain pFDR $\leq \alpha$. Then $n_* \sim (1/r) \ln Q_{\alpha,\pi}$ as $r \to 0+$.

**Remark.** We will refer to $r$ as the *signal-noise ratio* (SNR) of the multiple testing problem in Theorem 2.1.

Theorem 2.1 can be generalized to the case where the SNR follows a distribution. To specify how the SNR becomes increasingly small, we introduce a "scale" parameter $s > 0$ and parameterize the SNR distribution as $G_s(r) = G(sr)$, where $G$ is a fixed distribution.

**Corollary 2.1.** *Suppose that when $H_i : \mu_i = 0$ is false, $r_i = \mu_i/\sigma_i$ is a random sample from $G(sr)$, where $G(r)$ is a distribution function with support on $(0, \infty)$ and is unknown to the data analyzer. Suppose there is $\lambda > 0$, such that $\int e^{\lambda r^2} G(dr) < \infty$. Let $L_G$ be the Laplace transform of $G$, i.e., $L_G(\lambda) = \int e^{\lambda r} G(dr)$. Then $n_* \sim (1/s) L_G^{-1}(Q_{\alpha,\pi})$ as $s \to 0$.*

## 2.2. Preliminaries

Recall that, for the $t$ statistic (2.1), if $\mu = 0$, then $T \sim t_n$, the $t$ distribution with $n$ degrees of freedom (dfs). On the other hand, if $\mu > 0$, then $T \sim t_{n,\delta}$, the noncentral $t$ distribution with $n$ dfs and (noncentrality) parameter $\delta = \sqrt{n+1}\mu/\sigma$, with density

$$t_{n,\delta}(x) = \frac{n^{n/2}}{\sqrt{\pi}\,\Gamma(n/2)} \frac{e^{-\delta^2/2}}{(n+x^2)^{(n+1)/2}} \\ \times \sum_{k=0}^{\infty} \Gamma\left(\frac{n+k+1}{2}\right) \frac{(\delta x)^k}{k!} \left(\frac{2}{n+x^2}\right)^{k/2}.$$

Apparently $t_{n,0}(x) = t_n(x)$. Denote

$$a_{n,k} = \Gamma\left(\frac{n+k+1}{2}\right) \Big/ \Gamma\left(\frac{n+1}{2}\right).$$

Then

$$\frac{t_{n,\delta}(x)}{t_n(x)} = e^{-\delta^2/2} \sum_{k=0}^{\infty} \frac{a_{n,k}(\delta x)^k}{k!} \left(\frac{2}{n+x^2}\right)^{k/2}. \tag{2.2}$$

It can be shown that $t_{n,\delta}(x)/t_n(x)$ is strictly increasing in $x$ and

$$\sup_x \frac{t_{n,\delta}(x)}{t_n(x)} = \lim_{x \to \infty} \frac{t_{n,\delta}(x)}{t_n(x)} = e^{-\delta^2/2} \sum_{k=0}^{\infty} \frac{a_{n,k}(\sqrt{2}\,\delta)^k}{k!} < \infty \tag{2.3}$$

(cf. (3)). Since the supremum of likelihood ratio only depends on $n$ and $r = \mu/\sigma$, it will be denoted by $L(n,r)$ henceforth.

## 2.3. Proofs of the main results

We need two lemmas. They will be proved in the Appendix. The proofs of the main results are rather straightforward. The proofs are given in order illustrate the basic argument, which is used for the other results of the article as well.

**Lemma 2.1.** *1) For any fixed $n$, $L(n,r) \to 1$, as $r \to 0$. 2) Given $a \geq 0$, if $(n,r) \to (\infty, 0)$ such that $nr \to a$, then $L(n,r) \to e^a$. 3) If $(n,r) \to (\infty, 0)$ with $nr \to \infty$, then $L(n,r) \to \infty$.*

**Lemma 2.2.** *Under the same conditions as in Corollary 2.1, as $(n,s) \to (\infty, 0)$ such that $ns \to a \geq 0$, $\int L(n, sr) \, G(dr) \to L_G(a)$.*

*Proof of Theorem 2.1.* By (1.2), in order to get pFDR $\leq \alpha$,

$$\frac{1-\pi}{1-\pi+\pi L(n,r)} \leq \alpha, \quad \text{or} \quad L(n,r) \geq Q_{\alpha,\pi}.$$

Let $n_*$ be the minimum value of $n$ in order for the inequality to hold. Then by Lemma 2.1, as $r = \mu/\sigma \to 0$, $n_* r \to \ln Q_{\alpha,\pi}$, implying Theorem 2.1. □

*Proof of Corollary 2.1.* Following the argument for (1.2), it is seen that under the conditions of the corollary, the minimum attainable pFDR is

$$\alpha_* = \frac{1-\pi}{1-\pi+\pi \int L(n, sr) \, G(dr)}.$$

Then the corollary follows from a similar argument for Theorem (2.1). □

## 3. Multiple $F$-tests for linear regression with errors being normally distributed

### 3.1. Main results

Suppose we wish to test $H_i : \boldsymbol{\beta}_i = 0$ simultaneously for a large number of joint distributions of $Y$ and $\boldsymbol{X}$, such that under each distribution, $Y = \boldsymbol{\beta}_i^T \boldsymbol{X} + \epsilon_i$, where $\boldsymbol{\beta}_i \in \mathbb{R}^p$ are vectors of linear coefficients and $\epsilon_i \sim N(0, \sigma_i)$ are independent of $\boldsymbol{X}$. Suppose neither $\sigma_i$ or any possible relationships among $\sigma_i$ are known. Under this condition, consider the following tests based on a fixed design. Let $\boldsymbol{X}_k$, $k \geq 1$, be fixed vectors of covariates. Let $n + p$ be the sample size per null. For each $i$, let $(Y_{i1}, \boldsymbol{X}_1), \ldots, (Y_{i,n+p}, \boldsymbol{X}_{n+p})$ be an independent sample from $Y = \boldsymbol{\beta}_i^T \boldsymbol{X} + \epsilon$. Assume that the samples for different $H_i$ are independent of each other.

Suppose that, unknown to the data analyzer, for all the false nulls $H_i$,

$$\frac{(\boldsymbol{\beta}_i^T \boldsymbol{X}_1)^2 + \cdots + (\boldsymbol{\beta}_i^T \boldsymbol{X}_k)^2}{k \sigma_i^2} \leq \delta^2, \quad k = 1, 2, \ldots, \tag{3.1}$$

where $\delta > 0$. This situation arises when all $\boldsymbol{X}_k$ are within a bounded domain, either because only regression within the domain is of interest, or because only covariates within the domain are observable or experimentally controllable.

Note that $n$ is not the sample size per null. Instead, it is the difference between the sample size per null and the number of covariates in each regression equation. Given $\alpha \in (0,1)$, let

$$n_* = \inf\{n : \text{pFDR} \le \alpha \text{ for } F \text{ tests on } H_i \text{ under the constraint } (3.1)\}.$$

It can be seen that $n_*$ is attained when equality holds in (3.1) for all the false nulls. The asymptotics of $n_*$ will be considered for 3 cases: 1) $\delta \to 0$ while $p$ is fixed, 2) $\delta \to 0$ and $p \to \infty$, and 3) $p \to \infty$ while $\delta$ is fixed. The case $\delta \to 0$ is relevant when the regression effects are weak, and the case $p \to \infty$ is relevant when a large number of covariates are incorporated.

**Theorem 3.1.** *Under the random effects model* (1.1) *and the above setup of multiple $F$ tests, the following statements hold.*

a) *If $\delta \to 0$ while $p$ is fixed, then*

$$n_* \sim (1/\delta) M_p^{-1}(Q_{\alpha,\pi}), \quad \text{with} \quad M_p(t) := \sum_{k=0}^{\infty} \frac{\Gamma(p/2)(t^2/4)^k}{k!\Gamma(k+p/2)}.$$

b) *If $\delta \to 0$ and $p \to \infty$,*

$$n_* \sim \begin{cases} (1/\delta)\sqrt{2p \ln Q_{\alpha,\pi}} & \text{if } \delta^2 p \to 0, \\ (2/\delta^2) \ln Q_{\alpha,\pi} & \text{if } \delta^2 p \to \infty, \\ \dfrac{(4/\delta^2) \ln Q_{\alpha,\pi}}{1 + \sqrt{1 + 8 \ln Q_{\alpha,\pi}/L}} & \text{if } \delta^2 p \to L > 0. \end{cases}$$

c) *Finally, if $\delta > 0$ is fixed while $p \to \infty$, then*

$$n_* \to \left\lceil \frac{2 \ln Q_{\alpha,\pi}}{\ln(1+\delta^2)} \right\rceil.$$

### 3.2. Preliminaries and proofs

Given data $(Y_1, \boldsymbol{X}_1), \ldots, (Y_{n+p}, \boldsymbol{X}_{n+p})$, such that $Y_i = \boldsymbol{\beta}^T \boldsymbol{X}_i + \epsilon_i$, where $\boldsymbol{X}_i$ are fixed and $\epsilon_i$ are iid $\sim N(0, \sigma)$, if $\boldsymbol{\beta} = 0$, then the $F$ statistic of $(Y_i, \boldsymbol{X}_i)$ follows the $F$ distribution with $(p, n)$ dfs. On the other hand, if $\boldsymbol{\beta} \ne 0$, the $F$ statistic follows the noncentral $F$ distribution with $(p, n)$ dfs and (noncentrality) parameter $\Delta$, where

$$\Delta = \frac{(\boldsymbol{\beta}_i^T \boldsymbol{X}_1)^2 + \cdots + (\boldsymbol{\beta}_i^T \boldsymbol{X}_{n+p})^2}{\sigma_i^2}.$$

The density of the noncentral $F$ distribution is

$$f_{p,n,\Delta}(x) = e^{-\Delta/2}\theta^{p/2}x^{p/2-1}(1+\theta x)^{-(p+n)/2}$$
$$\times \sum_{k=0}^{\infty} \frac{(\Delta/2)^k}{k!\, B(p/2+k,\ n/2)} \left(\frac{\theta x}{1+\theta x}\right)^k, \qquad x \geq 0,$$

where $\theta = p/n$, and $B(a,b) = \Gamma(a)\Gamma(b)/\Gamma(a+b)$ is the Beta function. Note $f_{p,n,0}(x) = f_{p,n}(x)$, the density of the usual $F$ distribution with $(p,n)$ dfs.

Denote

$$b_{p,n,k} = \frac{B(p/2,\ n/2)}{B(p/2+k,\ n/2)} = \prod_{j=0}^{k-1}\left(\frac{n+p+2j}{p+2j}\right).$$

Then for $x \geq 0$,

$$\frac{f_{p,n,\Delta}(x)}{f_{p,n}(x)} = e^{-\Delta/2} \sum_{k=0}^{\infty} \frac{b_{p,n,k}(\Delta/2)^k}{k!} \left(\frac{\theta x}{1+\theta x}\right)^k, \qquad (3.2)$$

which is strictly increasing, and

$$\sup_{x>0} \frac{f_{p,n,\Delta}(x)}{f_{p,n}(x)} = \lim_{x\to\infty} \frac{f_{p,n,\Delta}(x)}{f_{p,n}(x)} = e^{-\Delta/2} \sum_{k=0}^{\infty} \frac{b_{p,n,k}(\Delta/2)^k}{k!} < \infty. \qquad (3.3)$$

First, it is easy to see that the following statement is true.

**Lemma 3.1.** *The expression in* (3.3) *is strictly increasing in* $\Delta > 0$.

It follows that, under the constraint (3.1), the supremum of the likelihood ratio is attained when $\Delta = (n+p)\delta^2$ and is equal to

$$K(p,n,\delta) = e^{-(n+p)\delta^2/2} \sum_{k=0}^{\infty} \frac{b_{p,n,k}[(n+p)\delta^2/2]^k}{k!}.$$

Therefore, under the random effects model (1.1), pFDR $\leq \alpha$ is equivalent to $K(p,n,\delta) \geq Q_{\alpha,\pi}$. Theorem 3.1 then follows from the lemmas below and an argument as to that of Theorem 2.1. The proof of Theorem 3.1 is omitted for brevity. The proofs of the lemmas are given in the Appendix.

**Lemma 3.2.** *Fix* $p \geq 1$. *If* $\delta \to 0$ *and* $n = n(\delta)$ *such that* $n\delta \to a \in [0,\infty)$, *then* $K(p,n,\delta) \to M_p(a)$. *If* $n\delta \to \infty$, *then* $K(p,n\delta) \to \infty$.

**Lemma 3.3.** *Let* $\delta \to 0$ *and* $p \to \infty$. *If* $n = n(\delta,p)$ *such that*

$$\frac{n(n+p)\delta^2}{2p} \to a \geq 0, \qquad (3.4)$$

*then* $K(p,n,\delta) \to e^a$. *In particular, given* $a > 0$, (3.4) *holds if*

$$n \sim \begin{cases} (1/\delta)\sqrt{2pa} & \text{if } \delta^2 p \to 0, \\ 2a/\delta^2 & \text{if } \delta^2 p \to \infty, \\ \dfrac{4a/\delta^2}{1+\sqrt{1+8a/L}} & \text{if } \delta^2 p \to L > 0. \end{cases}$$

**Lemma 3.4.** *Fix $\delta > 0$. Then for any $n \geq 1$, $K(n, p, \delta) \to (1 + \delta^2)^{n/2}$ as $p \to \infty$.*

## 4. Multiple $t$-tests: a general case

### 4.1. Setup

Suppose we wish to conduct hypothesis tests for a large number of distributions $F_i$ in order to identify those with nonzero mean $\mu_i$. The tests will be based on random samples from $F_i$. Assume that no information on the forms of $F_i$ or their relationships is available. As a result, samples from different $F_i$ cannot be combined to improve the inference. As in the case of testing mean values for normal distributions, to test $H_i : \mu_i = 0$ simultaneously, an appropriate approach is to use the $t$ statistics $T_i = \sqrt{n}\hat{\mu}_i/\hat{\sigma}_i$, where both $\hat{\mu}_i$ and $\hat{\sigma}_i^2$ are derived solely from the sample from $F_i$, and $n$ is the number of observations used to get $\hat{\mu}_i$.

Again, the goal is to find the minimum sample size per null in order to attain a given pFDR level, in particular when $F_i$ under false $H_i$ only have small differences from those under true $H_i$. The results will also answer the following question: are normal or $t$ approximations appropriate for the $t$ statistics in determining the minimum sample size per null?

We only consider the case where $\mu_i$ is either 0 or $\mu_0 \neq 0$, where $\mu_0$ is a constant. In order to make the analysis tractable, the problem needs to be formulated carefully. First, unlike the case of normal distributions, in general, if $\hat{\mu}_i$ and $\hat{\sigma}_i^2$ are the mean and variance of the same random sample, they are dependent and $\hat{\sigma}_i^2$ cannot be expressed as the sum of iid random variables. As seen below, the analysis on the minimum sample size per null requires detailed asymptotics of the $t$ statistics, in particular, the so called exact LDP (2, 5). For Studentized statistics, there are LDP techniques available (17). However, currently, *exact* LDP techniques cannot handle complex statistical dependency very well. To get around this technical difficulty, we consider the following $t$ statistics. Suppose the samples from different $F_i$ are independent of each other, and contain the same number of iid observations. Divide the sample from $F_i$ into two parts, $\{X_{i1}, \ldots, X_{in}\}$ and $\{Y_{i1}, Y_{i2}, \ldots, Y_{i,2m}\}$. Let

$$T_i = \frac{\sqrt{n}\hat{\mu}_i}{\hat{\sigma}_i}, \quad \text{with} \quad \hat{\mu}_i = \frac{1}{n}\sum_{k=1}^{n} X_{ik}, \quad \hat{\sigma}_i^2 = \frac{1}{2m}\sum_{k=1}^{n}(Y_{i,2k-1} - Y_{i,2k})^2.$$

Then $\hat{\mu}_i$ and $\hat{\sigma}_i^2$ are independent, and $\hat{\sigma}_i^2$ is the sum of iid random variables.

Second, the minimum attainable pFDR depends on the supremum of the ratio of the actual density of $T_i$ and its theoretical density under $H_i$. In general, neither one is tractable analytically. To deal with this difficulty, observe that in the case of normal distributions, the supremum of the ratio equals

$$\frac{P(T \geq t \mid \mu = \mu_0 > 0)}{P(T \geq t \mid \mu = 0)}, \qquad \text{as } t \to \infty.$$

We therefore consider the pFDR under the rule that $H_i$ is rejected if and only if $T_i > x$, where $x > 0$ is a critical value. In order to identify false nulls as $\mu_0 \to 0$, $x$ must increase, otherwise $P(T \geq x \,|\, \mu = \mu_0)/P(T \geq x \,|\, \mu = 0) \to 1$, giving pFDR $\to 1$. The question is how fast $x$ should increase.

Recall Section 2. Some analysis on (2.2) and (2.3) shows that, for normal distributions, the supremum of the likelihood ratio can be obtained asymptotically by letting $x = c_n\sqrt{n}$, where $c_n > 0$ is an arbitrary sequence converging to $\infty$; specifically, given $a > 0$, as $r \downarrow 0$ and $n \sim a/r$,

$$\frac{P(T > c_n\sqrt{n} \,|\, \mu/\sigma = r)/P(T > c_n\sqrt{n} \,|\, \mu = 0)}{\sup_x t_{n,r\sqrt{n}}(x)/t_n(x)} \to 1.$$

If, instead, $x$ increases in the same order as $\sqrt{n}$ or more slowly, the above limit is strictly less than 1. Based on this observation, for the general case, we set $x = c_n\sqrt{n}$, with $c_n \to \infty$. In general, there is no guarantee that using $c_n$ growing at a specific rate can always yield convergence. Thus, we require that $c_n$ grow slowly.

Under the setup, suppose that, unknown to the data analyzer, when $H_i : \mu_i = 0$ is true, $F_i(x) = F(s_i x)$, and when $H_i$ is false, $F_i(x) = F(s_i x - d)$, where $s_i > 0$ and $d > 0$, and $F$ is an unknown distribution such that

$$F \text{ has a density } f, \quad EX = 0, \quad \sigma^2 := EX^2 < \infty, \quad \text{for } X \sim F, \tag{4.1}$$

The sample from $F_i$ consists of $(X_{ij} - d)/s_i$, $1 \leq j \leq n$, and $(Y_{ik} - d)/s_i$, $1 \leq k \leq 2m$, with $X_{ij}, Y_{ik}$ iid $\sim F$. Then the $t$ statistic for $H_i$ is

$$T_i = \begin{cases} \sqrt{n}\bar{X}_{in}/S_{in} & \text{if } H_i \text{ is true,} \\ \sqrt{n}(\bar{X}_{in} + d)/S_{in} & \text{if } H_i \text{ is false,} \end{cases}$$

$$\text{where} \quad \bar{X}_{in} = \frac{X_{i1} + \ldots + X_{in}}{n}, \quad S_{im}^2 = \frac{1}{2m}\sum_{k=1}^{m}(Y_{i,2k-1} - Y_{i,2k})^2.$$

Let $N = n + m$ and $z_N = c_n$. Then $H_i$ is rejected if and only if $T_i \geq z_N\sqrt{n}$. Under the random effects model (1.1), the minimum attainable pFDR is

$$\alpha_* = (1 - \pi)\left[1 - \pi + \pi\frac{P(\bar{X}_n + d \geq z_N S_m)}{P(\bar{X}_n \geq z_N S_m)}\right]^{-1}, \tag{4.2}$$

where $\bar{X}_n = \sum_{k=1}^{n} X_k/n$, and $S_m = \sum_{k=1}^{m}(Y_{2k-1} - Y_{2k})^2/(2m)$, with $X_i, Y_j$ iid $\sim F$. The question now is the following:

- Given $\alpha \in (0, 1)$, as $d \to 0$, how should $N$ increase so that $\alpha_* \leq \alpha$?

### 4.2. Main results

By the Law of Large Numbers, as $n \to \infty$ and $m \to \infty$, $\bar{X}_n \to 0$ and $S_m \to \sigma$ w.p. 1. On the other hand, by our selection, $z_N \to \infty$. In order to analyze (4.2) as $d \to 0$, we shall rely on exact LDP, which depends on the properties of the cumulant generating functions

$$\Lambda(t) = \ln E e^{tX}, \quad \Psi(t) = \ln E\left[\exp\frac{t(X-Y)^2}{2}\right], \quad X, Y \text{ iid } \sim F. \quad (4.3)$$

The density of $X - Y$ is $g(t) = \int f(x)f(x+t)\,dx$. It is easy to see that $g(t) = g(-t)$ for $t > 0$. Recall that a function $\zeta$ is said to be slowly varying at $\infty$, if for all $t > 0$, $\lim_{x \to \infty} \zeta(tx)/\zeta(x) = 1$.

**Theorem 4.1.** *Suppose the following two conditions are satisfied.*

a) $0 \in \mathcal{D}_\Lambda^o$ *and* $\Lambda(t) \to \infty$ *as* $t \uparrow \sup \mathcal{D}_\Lambda$, *where* $\mathcal{D}_\Lambda = \{t : \Lambda(t) < \infty\}$.

b) *The density function $g$ is continuous and bounded on $(\epsilon, \infty)$ for any $\epsilon > 0$, and there exist a constant $\lambda > -1$ and a function $\zeta(z) \geq 0$ which is increasing in $z \geq 0$ and slowly varying at $\infty$, such that*

$$\lim_{x \downarrow 0} \frac{g(x)}{x^\lambda \zeta(1/x)} = C \in (0, \infty). \quad (4.4)$$

*Fix $\alpha \in (0, 1)$. Let $N_*$ be the minimum value for $N = m + n$ in order to attain $\alpha_* \leq \alpha$, where $\alpha_*$ is as in (4.2). Then, under the constraints 1) $m$ and $n$ grow in proportion to each other such that $m/N \to \rho \in (0, 1)$ as $m, n \to \infty$ and 2) $z_N \to \infty$ slowly enough, one gets*

$$N_* \sim \frac{1}{d} \times \frac{\ln Q_{\alpha,\pi}}{(1-\rho)t_0}, \qquad \text{as } d \to 0+, \quad (4.5)$$

*where $t_0 > 0$ is the unique positive solution to*

$$t\Lambda'(t) = \frac{(1+\lambda)\rho}{1-\rho}. \quad (4.6)$$

**Remark.** (1) By (4.5) and (4.6), $N_*$ depends on the moments of $F$ of all orders. Thus, $t$ or normal approximations of the distribution of $T$ in general are not suitable in determining $N_*$ in order to attain a target pFDR level.

(2) If $z_N \to \infty$ slowly enough such that (4.5) holds, then for any $z'_N \to \infty$ more slowly, (4.5) holds as well. Presumable, there is an upper bound for the growth rate of $z_N$ in order for (4.5) to hold. However, it is not available with the technique employed by this work.

(3) We define $N$ as $n + m$ instead of $n + 2m$ because in the estimator $S_m$, each pair of observations only generate one independent summand. The sum $n+m$ can be thought of as the number of degrees of freedom that are effectively utilized by $T$.

Following the proof for the case of normal distributions, Theorem 4.1 is a consequence of the following result.

**Proposition 4.1.** *Let $T > 0$. Under the same conditions as in Theorem 4.1, suppose $d = d_N \to 0$, such that $d_N N \to T > 0$. Then*

$$\frac{P\left(\bar{X}_n + d_N \geq z_N S_m\right)}{P\left(\bar{X}_n \geq z_N S_m\right)} \to e^{(1-\rho)Tt_0}. \tag{4.7}$$

Indeed, by display (4.2) and Proposition 4.1, if $dN \to T \geq 0$, then the minimum attainable pFDR has convergence

$$\alpha_* \to \frac{1-\pi}{1-\pi + \pi e^{(1-\rho)Tt_0}}. \tag{4.8}$$

In order to attain pFDR $\leq \alpha$, there must be $\alpha_* \leq \alpha$, leading to (4.5). The proof of Proposition 4.1 is given in the Appendix A3.

### 4.3. Examples

**Example 4.1** (Normal distribution)**.** Under the setup in Section 4.1, let $F = N(0, \sigma)$ in (4.1). By $\Lambda(t) = \ln E(e^{tX}) = \sigma^2 t^2/2$, condition a) of Theorem 4.1 is satisfied. For $X$, $Y$ iid $\sim F$, $X - Y \sim N(0, \sqrt{2}\sigma)$. Therefore, (4.4) is satisfied with $\lambda = 0$ and $\zeta(x) \equiv 1$. The solution to (4.6) is $t_0 = (1/\sigma)\sqrt{\rho/(1-\rho)}$. Then by Theorem 4.1,

$$N_* \sim \frac{\sigma}{d} \times \frac{\ln Q_{\alpha, \pi}}{\sqrt{\rho(1-\rho)}}, \qquad \text{as } d \to 0+. \tag{4.9}$$

To see the connection to Theorem 2.1, observe $\bar{X}_n = \sigma Z/\sqrt{n}$ and $S_m = \sigma W_m/\sqrt{m}$, where $Z \sim N(0,1)$ and $W_m^2 \sim \chi_m^2$ are independent. Since $z_N \uparrow \infty$ slowly, so is $a_m := \sqrt{n/m}\, z_N$. Let $r_m = (d/\sigma)\sqrt{n/(m+1)}$. Then

$$\frac{P(\bar{X}_n + d \geq z_N S_m)}{P(\bar{X}_n \geq z_N S_m)} = \frac{P(Z + \sqrt{m+1}\, r_m \geq a_m W_m)}{P(Z \geq a_m W_m)}$$
$$= \frac{1 - T_{m, \sqrt{m+1}\, r_m}(a_m)}{1 - T_m(a_N)},$$

where $T_{m,\delta}$ denotes the cumulative distribution function (cdf) of the noncentral $t$ distribution with $m$ dfs and parameter $\delta$, and $T_m$ the cdf of the $t$ distribution with $m$ dfs. Comparing the ratio in (2.2) and the above ratio, it is seen that the difference between the two is that probabilities densities in (2.2) are replaced with tail probabilities. Since $r_m = (d/\sigma)\sqrt{n/(m+1)} \sim (d/\sigma)\sqrt{(1-\rho)/\rho}$, by Theorem 2.1, in order to attain pFDR $\leq \alpha$ based on (2.2), the minimum value $m_*$ for $m$ satisfies $m_* \sim (\sigma/d)\sqrt{\rho/(1-\rho)} \ln Q_{\alpha, \pi}$. Since $m_*/N_* \to \rho$, the asymptotic of $N_*$ given by Theorem 2.1 is identical to that given by Theorem 4.1.

**Example 4.2** (Uniform distributions)**.** Under the setup in Section 4.1, let $F = U(-\frac{1}{2}, \frac{1}{2})$ in (4.1). Then for $t > 0$,

$$\Lambda(t) = -\frac{t}{2} + \ln(e^t - 1) - \ln t, \quad t\Lambda'(t) = \frac{t}{2 \tanh t} - 1.$$

and for $t < 0$, $\Lambda(t) = \Lambda(-t)$. Thus condition a) in Theorem 4.1 is satisfied. It is easy to see that condition b) is satisfied as well, with $\lambda = 0$ and $\zeta(x) \equiv 1$ in (4.4). Then by (4.5),

$$N_* \sim \frac{1}{d} \times \frac{\ln Q_{\alpha,\pi}}{2\tanh t_0}, \quad \text{with} \quad t_0 > 0 \quad \text{solving} \quad \frac{t}{\tanh t} = \frac{2}{1-\rho}. \qquad (4.10)$$

**Example 4.3** (Gamma distribution). Under the setup in Section 4.1, let $F$ be the distribution of $\xi - \alpha\beta$, where $\xi \sim \text{gamma}(\alpha, \beta)$ with density $\beta^{-\alpha} x^{\alpha-1} e^{-x/\beta}/\Gamma(\alpha)$. For $0 < t < 1/\beta$,

$$\Lambda(t) = \ln E[e^{t(\xi-\alpha\beta)}] = -\alpha\ln(1-\beta t) - \alpha\beta t, \quad t\Lambda'(t) = \frac{\alpha\beta^2 t^2}{1-\beta t}.$$

Therefore, condition a) in Theorem 4.1 is satisfied. Because the value of $\lambda$ in (4.4) is invariant to scaling, in order to verify condition b), without loss of generality, let $\beta = 1$. For $x > 0$, the density of $X - Y$ is then $g(x) = e^{-x} k(x)/\Gamma(\alpha)^2$, where $k(x) = \int_0^\infty u^{\alpha-1}(u+x)^{\alpha-1} e^{-2u}\, du$. It suffices to consider the behavior of $k(x)$ as $x \downarrow 0$. We need to analyze 3 cases.

**Case 1:** $\alpha > 1/2$ As $x \downarrow 0$, $k(x) \to \int_0^\infty u^{2\alpha-1} e^{-2u}\, du < \infty$. Therefore, (4.4) holds with $\lambda = 0$ and $\zeta \equiv 1$.

**Case 2:** $\alpha = 1/2$ As $x \downarrow 0$, $k(x) \to \infty$. We show that (4.4) still holds with $\lambda = 0$, but $\zeta(z) = \ln z$. To establish this, for any $\epsilon > 0$, let $k_\epsilon(x) = \int_0^\epsilon u^{-1/2}(u+x)^{-1/2}\, du$. Then

$$1 \leq \varliminf_{x\downarrow 0} \frac{k(x)}{k_\epsilon(x)} \leq \varlimsup_{x\downarrow 0} \frac{k(x)}{k_\epsilon(x)} \leq e^{2\epsilon}.$$

By variable substitution $u = xv^2$,

$$k_\epsilon(x) = 2 \int_0^{\sqrt{\epsilon/x}} \frac{dt}{\sqrt{t^2+1}} = (1+o(1))\ln(1/x), \qquad \text{as } x \downarrow 0.$$

As a result,

$$1 \leq \varliminf_{x\downarrow 0} \frac{k(x)}{\ln(1/x)} \leq \varlimsup_{x\downarrow 0} \frac{k(x)}{\ln(1/x)} \leq e^{2\epsilon}$$

Since $\epsilon$ is arbitrary, (4.4) is satisfied with $\lambda = 0$ and $\zeta(z) = \ln z$.

**Case 3:** $\alpha < 1/2$ As $x \downarrow 0$, $k(x) \to \infty$. Similar to the case $\alpha = -1/2$, it suffices to consider the behavior of $k_\epsilon(x) = \int_0^\epsilon u^{\alpha-1}(u+x)^{\alpha-1}\, du$ as $x \downarrow 0$, where $\epsilon > 0$ is arbitrary. By variable substitution $u = tx$,

$$k_\epsilon(x) = t^{2\alpha-1} \int_0^{\epsilon/x} t^{\alpha-1}(t+1)^{\alpha-1}\, dt = (1+o(1))C_\alpha t^{2\alpha-1}, \qquad \text{as } x \downarrow 0,$$

where $C_\alpha = \int_0^\infty t^{\alpha-1}(t+1)^{\alpha-1}\,dt < \infty$. Therefore, (4.4) is satisfied with $\lambda = 2\alpha - 1$ and $\zeta(z) \equiv 1$.

From the above analysis and (4.4), $N_* \sim (1/d)(\ln Q_{\alpha,\pi})/[(1-\rho)t_0]$, where

$$t_0 = \frac{\sqrt{\gamma^2 + 2\gamma} - \gamma}{\beta}, \quad \text{with} \quad \gamma = \frac{1}{1 \vee (2\alpha)} \times \frac{\rho}{1-\rho}. \tag{4.11}$$

### 4.4. Optimal split of sample

For the $t$ statistics considered so far, $m/N$ is the fraction of degrees of freedom allocated for the estimation of variance. By (4.5), the asymptotic of $N_*$ depends on the fraction in a nontrivial way. It is of interest to optimize the fraction in order to minimize $N_*$. Asymptotically, this is equivalent to maximizing $(1-\rho)t_0$ as a function of $\rho$, with $t_0 = t_0(\rho) > 0$ as in (4.6).

**Example** 4.1 (Continued)

By (4.9), it is apparent that the optimal value of $\rho$ is $1/2$. In other words, in order to minimize $N_*$, there should be equal number of degrees of freedom allocated for the estimation of mean and the estimation of variance for each normal distribution. In particular, if $m \equiv n-1$, then $\rho = 1/2$, and the resulting $t$ statistic has the same distribution as $\sqrt{n-1}Z/W_{n-1}$, where $Z \sim N(0,1)$ and $W_{n-1} \sim \chi_{n-1}$ are independent, which is the usual $t$ statistic of an iid sample of size $n$.

**Example** 4.2 (Continued)

By (4.10), the larger $\tanh t_0$ is, the smaller $N_*$ becomes. The function $\tanh t_0$ is strictly increasing in $t_0$, and $\tanh t_0 \to 1$ as $t_0 \to \infty$. By $\rho = 1 - 2\tanh t_0/t_0$, the closer $\rho$ is to 1, the smaller $N_*$.

**Example** 4.3 (Continued)

Denote $\theta = 1/[1 \vee (2\alpha)]$. By (4.11), we need to find $\rho$ to maximize

$$(1-\rho)\left[\sqrt{\gamma^2 + 2\gamma} - \gamma\right] = \sqrt{\theta^2\rho^2 + 2\theta\rho(1-\rho)} - \theta\rho.$$

By some calculation, the value of $\rho$ that maximizes the above quantity is

$$\rho_* = \frac{1}{2 + \sqrt{2\theta}} = \frac{1}{2 + \sqrt{2 \wedge (1/\alpha)}}.$$

For $0 < \alpha \leq 1/2$, the optimal fraction of degrees of freedom allocated for the estimation of the variance of gamma$(\alpha, \beta)$ tends to $1/(2 + \sqrt{2})$ as $d \to 0$. On the other hand, as $\alpha \to \infty$, the optimal fraction tends to $1/2$ as $d \to 0$, which is reasonable in light of Example 4.1. To see this, let $\beta = 1$. For integer valued $\alpha$ and $\xi \sim \text{gamma}(\alpha, 1)$, $\xi - \alpha$ can be regarded as the sum of $W_i - 1$, $i = 1, \ldots, \alpha$, with $W_i$ iid following gamma$(1, 1)$. Therefore, for $\alpha \gg 1$, $\xi - \alpha$ follows closely a normal distribution with mean 0. Thus by Example 4.1, the optimal value of $m/(n + m)$ is close to $1/2$.

## 5. Multiple tests based on likelihoods

### 5.1. Motivation

In many cases of multiple testing, only limited knowledge is available on the distributions from which data are sampled. The knowledge relevant to a null hypothesis is expressed by a statistic $M$ such that the null is rejected if and only if the observed value of $M$ is significantly different from 0. In general, as the distribution of $M$ is unknown, $M$ has to be Studentized so that its magnitude can be evaluated.

On the other hand, oftentimes, despite the complexity of the data distributions, it is reasonable to believe they have an underlying structure. Consider the scenario where all the data distributions belong to a parametric family $\{p_\theta\}$, such that the distribution under a true null is $p_0$, and the one under a false null is $p_{\theta_*}$ for some $\theta_* \neq 0$. A question of interest is: under this circumstance, what would be the optimal overall performance of the multiple tests? The question is in the same spirit as questions regarding estimation efficiency. However, it assumes that neither the existence of the parameterization nor its form is known to the data analyzer and all the machinery available is the test statistic $M$.

As before, we wish to find out the minimum sample size per null required for pFDR control, in particular, as the tests become increasingly harder in the sense that $\theta_* \to 0$. Our conjecture is that, asymptotically the minimum sample size per null is attained if $M$ "happens" to be $\partial[\ln p_0]/\partial\theta$. By "happens" we mean that the data analyzer is unaware of this peculiar nature of $M$ and uses its Studentized version for the tests. This conjecture is directly motivated by the fact that the MLE is efficient under regular conditions. Although a smaller minimum sample size per null could be possible if $M$ happens to be the MLE, due to Studentization, the improvement appears to diminish as $\theta \to 0$. Certainly, had the parameterization been known, the (original) MLE would be preferred. The goal here is not to establish any sort of superiority of Studentized MLE, but rather to search for the optimal overall performance of multiple tests, when we are aware that our knowledge about the data distributions is incomplete and beyond the test statistic, we have no other information.

The above conjecture is not yet proved or disproved. However, as a first step, we would like to obtain the asymptotics of the minimum sample size per null when Studentized $\partial[\ln p_0]/\partial\theta$ is used for multiple tests. We shall also provide some examples to support the conjecture.

### 5.2. Setup

Let $(\Omega, \mathcal{F})$ be a measurable space equipped with a $\sigma$-finite measure $\mu$. Let $\{p_\theta : \theta \in [0,1]\}$ be a parametric family of density functions on $(\Omega, \mathcal{F})$ with respect to $\mu$. Denote by $P_\theta$ the corresponding probability measure. Under the random effects model (1.1), each null $H_i$ is associated with a distribution $F_i$, such that when $H_i$ is true, $F_i = P_0$, and when $H_i$ is false, $F_i = P_\theta$, where $\theta > 0$ is a constant. Assume that each $H_i$ is tested based on an iid sample $\{\omega_{ij}\}$ from $F_i$, such that the samples for different $H_i$ are independent, and the sample size is the same for all $H_i$.

We need to assume some regularities for $p_\theta$. Denote

$$r_\theta(\omega) = \frac{p_\theta(\omega)}{p_0(\omega)}, \quad \ell_\theta(\omega) = \ln p_\theta(\omega), \qquad \omega \in \Omega. \tag{5.1}$$

**Condition 1** Under $P_0$, for almost every $\omega \in \Omega$, $p_0(\omega) > 0$ and $p_\theta(\omega)$ as a function of $\theta$ is in $C^2([0,1])$.

**Condition 2** The Fisher information at $\theta = 0$ is positive and finite, i.e. $0 < \|\dot\ell_0\|_{L^2(P_0)} < \infty$, where the "dot" notation denotes partial differentiation with respect of $\theta$.

**Condition 3** Under $P_0$, the second order derivative of $\ell_\theta(\omega)$ is uniformly bounded in the sense that $\sup_{\theta \in [0,1]} \|\ddot\ell_\theta(\omega)\|_{L^\infty(P_0)} < \infty$.

**Condition 4** For any $q > 0$, there is $\theta' = \theta'(q) > 0$, such that

$$E_0\left[\sup_{\theta \in [0, \theta']} (r_\theta(\omega)^q + r_\theta(\omega)^{-q})\right] < \infty. \tag{5.2}$$

**Remark.** By Condition 1, for any interval $I$ in $[0,1]$, the extrema of $r_\theta(\omega)$ over $\theta \in I$ are measurable. Thus the expectation in (5.2) is well defined.

For brevity, for $\theta \in [0,1]$ and $n \geq 1$, the $n$-fold product measure of $P_\theta$ is still denoted by $P_\theta$, and the expectation under the product measure by $E_\theta$. We shall denote by $\omega$, $\omega'$, $\omega_i$, $\omega_i'$ generic iid elements under a generic distribution on $(\Omega, \mathcal{F})$. Denote

$$X = \dot\ell_0(\omega), \ Y = \dot\ell_0(\omega'), \ X_i = \dot\ell_0(\omega_i), \ Y_i = \dot\ell_0(\omega_i'). \tag{5.3}$$

For $m, n \geq 1$, denote

$$S_m^2 = \frac{1}{m} \sum_{i=1}^m \frac{(Y_{2i-1} - Y_{2i})^2}{2}, \quad \bar X_n = \frac{X_1 + \cdots + X_n}{n}.$$

Since $\dot\ell_\theta(\omega) = \dot p_\theta(\omega)/p_\theta(\omega)$, from Conditions 1–4 and dominated convergence, it follows that $E_0 \dot\ell_0 = 0$ and

$$(E_\theta \dot\ell_0)'\big|_{\theta=0} = \frac{d}{d\theta} \int \dot\ell_0(\omega) p_\theta(\omega)\,\mu(d\omega)\bigg|_{\theta=0} = \int (\dot\ell_0(\omega))^2 p_0(\omega)\,\mu(d\omega) > 0.$$

As a result, for $\theta > 0$ close to 0, $E_\theta \dot\ell_0(\omega) > 0$. This justifies using the upper tail of $\sqrt{n}\bar{X}_n/S_m$ for testing. The multiple tests are such that

$$H_i \text{ is rejected} \iff \frac{\sqrt{n}\bar{X}_{in}}{S_{im}} \geq z_N\sqrt{n}, \qquad i \geq 1, \qquad (5.4)$$

where $\bar{X}_{in}$ and $S_{im}$ are computed the same way as $\bar{X}_n$ and $S_m$, except that they are derived from $\omega_{i1}, \ldots, \omega_{in}, \omega'_{i1}, \ldots, \omega'_{i,2m}$ iid $\sim F_i$, $N = n+m$, and $z_N \to \infty$ as $N \to \infty$. Then, under the random effects model, the minimum attainable pFDR is

$$\alpha_* = (1-\pi)\left[1 - \pi + \pi\frac{P_\theta\left(\bar{X}_n \geq z_N S_m\right)}{P_0\left(\bar{X}_n \geq z_N S_m\right)}\right]^{-1}. \qquad (5.5)$$

The question now is the following:

- Given $\alpha \in (0,1)$, as $\theta \downarrow 0$, how should $N$ increase so that $\alpha_* \leq \alpha$?

### 5.3. Main results

Denote the cumulant generating functions

$$\Lambda(t) = \ln E_0(e^{tX}), \qquad \Psi(t) = \ln E_0\left[\exp\frac{t(X-Y)^2}{2}\right]. \qquad (5.6)$$

Note that the expectation is taken under $P_0$.

**Theorem 5.1.** *Suppose $\{p_\theta : \theta \in [0,1]\}$ satisfies conditions 1–4 and the following conditions a)–d) are fulfilled..*

a) $0 \in \mathcal{D}_\Lambda^o$, where $\mathcal{D}_\Lambda = \{t : \Lambda(t) < \infty\}$.

b) Under $P_0$, $X$ has a density $f$ continuous almost everywhere on $\mathbb{R}$. Furthermore, either (i) $f$ is bounded or (ii) $f$ is symmetric and $\|X\|_{L^\infty(P_0)} < \infty$.

c) Under $P_0$, the density $g$ of $X - Y$ is continuous and bounded on $(\epsilon, \infty)$ for any $\epsilon > 0$, and there exist a constant $\lambda > -1$ and a function $\zeta(z) \geq 0$ increasing in $z \geq 0$ and slowly varying at $\infty$, such that

$$\lim_{u \downarrow 0} \frac{g(u)}{u^\lambda \zeta(1/u)} = C \in (0,\infty). \qquad (5.7)$$

d) There are $s > 0$ and $L > 0$, such that

$$E_0[e^{s|X+Y|} \mid X - Y = u] \leq Le^{L|u|}, \qquad \text{any } u \neq 0,\ g(u) > 0. \qquad (5.8)$$

*Fix $\alpha \in (0,1)$. Let $N_*$ be the minimum value of $N = n+m$ in order to attain $\alpha_* \leq \alpha$, where $\alpha_*$ is as in (5.5). Then, under the constraints 1) $m$ and $n$ grow* *in proportion to each other such that $m/N \to \rho \in (0,1)$ as $m, n \to \infty$ and 2) $z_N \to \infty$ slowly enough, one gets*

$$N_* \sim \frac{1}{d} \times \frac{\ln Q_{\alpha,\pi}}{(1-\rho)\Lambda'(t_0) + 2\rho K_f}, \qquad \text{as } d \to 0+. \tag{5.9}$$

*where $t_0$ is the unique positive solution to (4.6), and*

$$K_f = \begin{cases} \int z h_0(z)\,dz & \text{if } f \text{ is bounded, with } h_0 = f^2 / \int f^2, \\ 0 & \text{if } f \text{ is symmetric and } \|X\|_{L^\infty(P_0)} < \infty. \end{cases}$$

**Remark.** By symmetry, to verify (5.8), it is enough to only consider $u > 0$. Moreover, (5.8) holds if its left hand side is a bounded function of $u$.

Following the proofs of the previous results, Theorem 5.1 is a consequence of Proposition 5.1, which will be proved in Appendix A4.

**Proposition 5.1.** *Let $T > 0$. Under the same conditions as in Theorem 5.1, suppose $\theta = \theta_N \to 0$, such that $\theta_N N \to T$. Then*

$$\frac{P_{\theta_N}(\bar{X}_n \geq z_N S_m)}{P_0(\bar{X}_n \geq z_N S_m)} \to \exp\{(1-\rho)T\Lambda'(t_0) + 2\rho T K_f\}. \tag{5.10}$$

### 5.4. Examples

**Example 5.1** (Normal distributions). Under the setup in Section 5.2, suppose for $\theta \in [0,1]$, $P_\theta = N(\theta, \sigma)$, where $\sigma > 0$ is a fixed constant. Then $p_\theta(u) = \exp[-(u-\theta)^2/(2\sigma^2)]/\sqrt{2\pi\sigma^2}$, $u \in \mathbb{R}$, giving

$$r_\theta(u) = \exp\left(\frac{2\theta u - \theta^2}{\sigma^2}\right), \quad \ell_\theta(u) = -\frac{(u-\theta)^2}{2\sigma^2} - \frac{\ln(2\pi\sigma^2)}{2},$$

$$\dot{\ell}_\theta(u) = \frac{u-\theta}{\sigma^2}, \quad \ddot{\ell}_\theta(u) = -\frac{1}{\sigma^2}.$$

For $\omega \sim P_0$, $\dot{\ell}_0(\omega) = \omega/\sigma^2 \sim N(0, 1/\sigma)$. It is then not hard to see that Conditions 1–4 are satisfied. By the notations in (5.3), $X$, $Y$, $X_i$, $Y_i$ are iid $\sim N(0, 1/\sigma)$. Then $\Lambda(t) = t^2/(2\sigma^2)$ and condition a) of Theorem 5.1 is satisfied. It is easy to see that conditions b) and c) are satisfied with $\lambda = 0$ and $\zeta \equiv 1$ in (5.7). Since $X + Y$ and $X - Y$ are independent and the moment generating function of $|X + Y| \sim \sqrt{2}|X|$ is finite on the entire $\mathbb{R}$, (5.8) is satisfied as well. Therefore, (5.10) holds. Therefore, (5.9) holds for $N_*$.

To get the asymptotic in (5.9) explicitly, note that the density $f$ of $X$ is $p_0$. Then it is not hard to see $K_f = 0$. On the other hand, since $\Lambda'(t) = t/\sigma^2$, the solution $t_0 > 0$ to $t\Lambda'(t) = \sqrt{\rho/(1-\rho)}$ equals $\sigma\sqrt{\rho/(1-\rho)}$ and hence $\Lambda'(t_0) = (1/\sigma)\sqrt{\rho/(1-\rho)}$. Thus, $N_* \sim (\sigma/d)(\ln Q_{\alpha,\pi}/\sqrt{\rho(1-\rho)})$, which is identical to (4.9) for the $t$ tests.

**Example 5.2** (Cauchy distribution). Under the setup in Section 5.2, suppose for $\theta \in [0,1]$, $P_\theta$ is the Cauchy distribution centered at $\theta$ such that its density is $p_\theta(u) = \pi^{-1}[1 + (u-\theta)^2]^{-1}$, $u \in \mathbb{R}$. Then

$$r_\theta(u) = \frac{1+u^2}{1+(u-\theta)^2}, \quad \ell_\theta(u) = -\ln[1+(u-\theta)^2] - \ln \pi,$$

$$\dot\ell_\theta(u) = \frac{2(u-\theta)}{1+(u-\theta)^2}, \quad \dot\ell_0(u) = \frac{2u}{1+u^2}.$$

By the notations in (5.3), $X = 2\omega/(1+\omega^2)$, with $\omega \sim P_0$. Recall that $P_0$ is the distribution of $\tan(\xi/2)$ with $\xi \sim U(-\pi,\pi)$. Therefore, $X \sim \sin \xi$ and thus is bounded and has a symmetric distribution. It is clear that conditions a), b), and d) of Theorem 5.1 are satisfied. We show that condition c) is satisfied with $\lambda = 0$ and $\zeta(z) = \ln z$ in (5.7). The density $f$ of $X$ is $1/[\pi\sqrt{1-u^2}]$, $u \in [-1,1]$. Then $K_f = 0$ and the density of $X - Y$ is

$$g(u) = k(u)/\pi^2, \quad \text{with} \quad k(u) = \int_{-1}^{1-u} \frac{dt}{\sqrt{(1-t^2)[1-(t+u)^2]}}, \quad u \in (0,1).$$

Given $\epsilon \in (0, 1-u/2)$, write the integral as the sum of integrals over $[-1, -1+\epsilon]$, $[1-u-\epsilon, 1-u]$, and $[-1+\epsilon, 1-u-\epsilon]$. By variable substitution

$$k(u) = 2\int_0^\epsilon \frac{dt}{\sqrt{(2-t)(2-t-u)t(t+u)}} + \int_{-1+\epsilon}^{1-u-\epsilon} \frac{dt}{\sqrt{(1-t^2)[1-(t+u)^2]}}$$

$$\sim 2\int_0^\epsilon \frac{dt}{\sqrt{(2-t)(2-t-u)t(t+u)}}, \quad \text{as} \quad u \to 0.$$

Because $\epsilon > 0$ is arbitrary, it follows that $k(u) \sim k_1(u)$, where

$$k_1(u) = \int_0^\epsilon \frac{dt}{\sqrt{t(t+u)}} = 2\int_0^{\sqrt{\epsilon/u}} \frac{dx}{\sqrt{x^2+1}} \sim \ln(1/u),$$

with the second equality due to variable substitution $t = ux^2$. This shows that (5.7) holds with $\lambda = 0$ and $\zeta(z) = \ln z$. By (5.9), $N_* \sim (t_0/d) \times (\ln Q_{\alpha,\pi}/\rho)$, as $d \to 0$, where $t_0$ the positive solution to $t_0 \Lambda'(t_0) = \rho/(1-\rho)$, with $\Lambda(t) = \ln E[e^{t \sin \xi}]$.

**Remark.** Because the Cauchy distributions have infinite variance, $t$ tests cannot be used to test the nulls. The example shows that even in this case, Studentized $\ell_0(\omega)$ can still distinguish between true and false nulls.

**Example 5.3** (Gamma distribution). Under the setup in Section 5.2, suppose for $\theta \in [0,1]$, $P_\theta = \text{gamma}(1+\theta, 1)$, whose density is $p_\theta(u) = u^\theta e^{-u}/\Gamma(1+\theta)$, $u > 0$. Then

$$r_\theta(u) = \frac{u^\theta}{\Gamma(\theta+1)}, \quad \ell_\theta(u) = \theta \ln u - u - \ln \Gamma(\theta+1),$$

$$\dot\ell_\theta(u) = \ln u - \psi(\theta+1), \quad \ddot\ell_\theta(u) = -\psi'(\theta+1),$$

where $\psi(z) = \Gamma'(z)/\Gamma(z)$ is the digamma function. Let $c = \psi(1)$. By the notations in (5.3), $X$ and $Y$ are iid $\sim \ln\omega - c$, with $\omega \sim P_0$. It follows that $X$ has density $f(x) = e^{x+c}p_0(e^{x+c}) = e^{x+c}\exp(-e^{x+c})$, $x \in \mathbb{R}$, which is bounded and continuous, and hence conditions b) and c) of Theorem 5.1 are satisfied with $\lambda = 1$ and $\zeta(z) = 1$ in (5.7). Since

$$E_0[e^{tX}] = \int_{-\infty}^{\infty} e^{tx}e^{x+c}\exp(-e^{x+c})\,dx$$
$$= \int_0^{\infty} z^t e^c \exp(-e^c z)\,dz = \frac{\Gamma(t+1)}{e^{ct}} < \infty, \quad \text{any} \quad t > -1,$$

condition a) is satisfied. To verify d), the density of $X - Y$ at $u > 0$ is

$$g(u) = \int_{-\infty}^{\infty} e^{2c+u+2x}\exp\left[-(1+e^u)e^{c+x}\right]\,dx \qquad \text{(substitute } z = e^{c+x})$$
$$= e^u \int_0^{\infty} z\exp[-(1+e^u)z]\,dz = \frac{e^u}{(1+e^u)^2}.$$

Similarly, for $s > 0$,

$$k(s,u) := \int_{-\infty}^{\infty} e^{s(2x+u)}e^{2c+u+2x}\exp\left[-(1+e^u)e^{c+x}\right]\,dx$$
$$= e^{(1+s)u-2sc}\int_0^{\infty} z^{1+2s}\exp[-(1+e^u)z]\,dz = \frac{\Gamma(2+2s)e^{(1+s)u}}{e^{2sc}(1+e^u)^{2+2s}}.$$

As a result, for $s \leq 1/2$,

$$E_0[e^{s(X+Y)} \mid X - Y = u] = \frac{k(s,u)}{g(u)} = \Gamma(2+2s)\left[\frac{e^u}{e^{2c}(1+e^u)^2}\right]^s.$$

Likewise,

$$E_0[e^{-s(X+Y)} \mid X - Y = u] = \Gamma(2-2s)\left[\frac{e^u}{e^{2c}(1+e^u)^2}\right]^{-s}.$$

Since $e^{s|X+Y|} \leq e^{s(X+Y)} + e^{-s(X+Y)}$, it is not hard to see that we can choose $s = 1/2$ and $L > 0$ large enough, such that (5.8) holds.

By $\Lambda(t) = \ln\Gamma(t+1) - \psi(1)t$, $t_0 > 0$ is the solution to $t[\psi(t+1) - \psi(1)] = \rho/(1-\rho)$. By $\int f^2 = g(0) = 1/4$,

$$K_f = 4\int_{-\infty}^{\infty} zf(z)^2\,dz = 4\int_{-\infty}^{\infty} ze^{2z+2c}\exp(-2e^{z+c})\,dz,$$

which equals $\psi(2) - \ln 2 - \psi(1)$. By $\psi(z) = (\ln\Gamma(z))'$ and $\Gamma(z+1) = z\Gamma(z)$, $\psi(z+1) - \psi(z) = 1/z$. Therefore, $K_f = 1 - \ln 2$. So by (5.9), $N_* \sim (1/d) \times \ln Q_{\alpha,\pi}/[(1-\rho)\Lambda'(t_0) + 2\rho(1-\ln 2)]$.

## 6. Summary

Multiple testing is often used to identify subtle real signals (false nulls) from a large and relatively strong background of noise (true nulls). In order to have some assurance that there is a reasonable fraction of real signals among the signals "spotted" by a multiple testing procedure, it is useful to evaluate the pFDR of the procedure. Comparing to FDR control, pFDR control is more subtle and in general requires more data. In this article, we study the minimum number of observations per null in order to attain a target pFDR level and show that it depends on several factors: 1) the target pFDR control level, 2) the proportion of false nulls among the nulls being tested, 3) distributional properties of the data in addition to mean and variance, and 4) in the case of multiple $F$ tests, the number of covariates included in the nulls.

The results of the article indicate that, in determining how much data are needed for pFDR control, if there is little information about the data distributions, then it may be useful to estimate the cumulant generating functions of the distributions. Alternatively, if one has good evidence about the parametric form of the data distributions but has little information on the values of the parameters, then it may be necessary to determine the number of observations per null based on the cumulant functions as well. In either case, typically it is insufficient to only use the means and variances of the distributions.

The article only considers univariate test statistics, which allow detailed analysis of tail probabilities. It is possible to test each null by more than one statistic. How to determine the number of observations per null for multivariate test statistics is yet to be addressed.

## Appendix: Mathematical Proofs

### A1. Proofs for normal *t*-tests

*Proof of Lemma 2.1.* Part 1) is clear. To show 2), let $(n, r) \to (\infty, 0)$ such that $nr \to a \geq 0$. Since $\delta = \sqrt{n+1}\, r \to 0$, by (2.3), it suffices to show

$$\sum_{k=0}^{\infty} \frac{a_{n,k}(\sqrt{2}\,\delta)^k}{k!} \to e^a. \tag{A1.1}$$

By Stirling's formula, $\Gamma(x) = (z/e)^z \sqrt{2\pi/z}\,[1 + O(1/z)]$. Then for $n \gg 1$,

$$a_{n,k} \leq 2\left(\frac{n+k+1}{2e}\right)^{(n+k+1)/2} \left(\frac{n+1}{2e}\right)^{-(n+1)/2}$$

$$\leq 2\left(\frac{n+k+1}{2e}\right)^{k/2} \left(1 + \frac{k}{n+1}\right)^{(n+1)/2} = 2\left(\frac{n+k+1}{2}\right)^{k/2},$$

giving

$$\frac{a_{n,k}(\sqrt{2}\,\delta)^k}{k!} \leq \frac{2(\sqrt{2}\,\delta)^k}{k!}\left(\frac{n+k+1}{2}\right)^{k/2}$$
$$= \frac{2[(n+1)(n+1+k)r^2]^{k/2}}{k!} \leq \frac{3(nr+r)^k(1+k)^{k/2}}{k!}. \tag{A1.2}$$

The right hand side has a finite sum over $k$. By dominated convergence,

$$\lim_{\substack{(n,r)\to(\infty,0)\\ \text{s.t. } nr\to a}} L(n,r) = \sum_{k=0}^{\infty} \lim_{\substack{(n,r)\to(\infty,0)\\ \text{s.t. } nr\to a}} \frac{a_{n,k}(\sqrt{2}\,\delta)^k}{k!}$$
$$= \sum_{k=0}^{\infty} \lim_{\substack{(n,r)\to(\infty,0)\\ \text{s.t. } nr\to a}} \frac{[(n+1)(n+1+k)r^2]^{k/2}}{k!} = \sum_{k=0}^{\infty} \frac{a^k}{k!} = e^a.$$

This yields 2). To show 3), by similar argument, given $0 < c < 1$, for $n \gg 1$,

$$\frac{a_{n,k}(\sqrt{2}\,\delta)^k}{k!} \geq \frac{c(\sqrt{2}\,\delta)^k}{k!}\left(\frac{n+1}{2}\right)^{k/2} \geq \frac{c(nr)^k}{k!}$$

Therefore, as $nr \to \infty$, $L(n,r) \geq c e^{nr} \to \infty$. $\square$

**Proof of Lemma 2.2**

By Stirling's formula, there is a constant $C > 1$, such that $k^{k/2}/k! \leq C^k/\Gamma(k/2+1)$ for all $k \gg 1$. Fix $n_0$ so that $C^2 a^2/n_0 < \lambda$ and (A1.2) holds for all $n \geq n_0$. For $k \geq n_0(n_0+1)$, $1 + k/(n_0+1) \leq k/n_0$. Then applying (A1.2) with $\delta = \sqrt{n+1}\,sr$ yields

$$\frac{a_{n,k}(\sqrt{2}\delta)^k}{k!} \leq \frac{2[(n+1)(n+1+k)s^2 r^2]^{k/2}}{k!} \leq \frac{2(k/n_0)^{k/2}(nsr+sr)^k}{k!}$$
$$\leq \frac{2\,C^k}{\Gamma(k/2+1)}\left[\frac{(ns+s)^2 r^2}{n_0}\right]^{k/2} \leq \frac{2[b(s)r^2]^{k/2}}{\Gamma(\lfloor k/2\rfloor + 1)},$$

where $b(s) = C^2(ns+s)^2/n_0$. Let $\lambda_* \in (C^2 a^2/n_0, \lambda)$. By $\int e^{\lambda r^2} G(dr) < \infty$, $\int r^p e^{\lambda_* r^2} G(dr) < \infty$ for any $p \geq 0$. Let $(n,s) \to (\infty, 0)$ such that $ns \to a$. Then for $n \gg n_0$, $b(s) < \lambda_*$ and hence

$$\sum_{k=k_0}^{\infty} \frac{a_{n,k}(\sqrt{2(n+1)}st)^k}{k!} \leq \sum_{k=k_0}^{\infty} \frac{2(b(s)r^2)^{k/2}}{\Gamma(\lfloor k/2 \rfloor + 1)}$$

$$\leq 2(1 + \sqrt{\lambda_*}r) \sum_{k=\lfloor k_0/2 \rfloor}^{\infty} \frac{(\lambda_* r^2)^k}{k!}.$$

By the above inequality and dominated convergence,

$$\lim \int L(n, sr)\, G(dr) = \int L(n, sr)\, G(dr) = \int e^{ar}\, G(dr). \qquad \square$$

## A2. Proofs for $F$-tests

**Proof of Lemma 3.1**

It suffices to show $\phi'(t) > 0$ for $t > 0$, where

$$\phi(t) = e^{-t} \sum_{k=0}^{\infty} \frac{b_{p,n,k} t^k}{k!}.$$

This follows from $b_{p,n,k+1} > b_{p,n,k}$ and

$$\phi'(t) = -\phi(t) + e^{-t} \sum_{k=0}^{\infty} \frac{b_{p,n,k+1} t^k}{k!} = e^{-t} \sum_{k=0}^{\infty} \frac{[b_{p,n,k+1} - b_{p,n,k}] t^k}{k!} > 0. \qquad \square$$

Next, recall

$$K(p, n, \delta) = e^{-(n+p)\delta^2/2} \sum_{k=0}^{\infty} \prod_{j=0}^{k-1} \left( \frac{n+p+2j}{p+2j} \right) \times \frac{1}{k!} \left[ \frac{(n+p)\delta^2}{2} \right]^k.$$

*Proof of Lemma 3.2.* Suppose $\delta \to \infty$ and $n = n(\delta)$ such that $n\delta \to a \in [0, \infty)$. Since $(n+p+2j)/(p+2j) \leq n+p$, then

$$K(p, n, \delta) \leq \sum_{k=0}^{\infty} (n+p)^k \frac{1}{k!} \left[ \frac{(n+p)\delta^2}{2} \right]^k \leq e^{(n+p)^2 \delta^2/2}, \qquad (A2.1)$$

and by dominated converge,

$$\lim_{\delta \to \infty} K(p, n, \delta) = \lim_{\delta \to \infty} \sum_{k=0}^{\infty} \prod_{j=0}^{k-1} \left[ \frac{1 + 2j/(n+p)}{p+2j} \right] \times \frac{1}{k!} \left[ \frac{(n+p)^2 \delta^2}{2} \right]^k$$

$$= \sum_{k=0}^{\infty} \prod_{j=0}^{k-1} \left( \frac{1}{p+2j} \right) \times \frac{1}{k!} \left( \frac{a^2}{2} \right)^k = M_p(a).$$

Next suppose $\delta \to 0$ and $n\delta \to \infty$. Then one gets

$$K(p, n, \delta) \geq e^{-(n+p)\delta^2/2} \sum_{k=0}^{\infty} \prod_{j=0}^{k-1} \left(\frac{n}{p+2j}\right) \times \frac{1}{k!}\left(\frac{n\delta^2}{2}\right)^k$$

$$\geq e^{-(n+p)\delta^2/2} \sum_{k=0}^{\infty} \prod_{j=0}^{k-1} \left(\frac{1}{1+2j}\right) \times \left(\frac{n}{p}\right)^k \times \frac{1}{k!}\left(\frac{n\delta^2}{2}\right)^k$$

$$= e^{-(n+p)\delta^2/2} \sum_{k=0}^{\infty} \frac{1}{(2k)!}\left(\frac{n^2\delta^2}{p}\right)^k = \frac{e^{-(n+p)\delta^2/2}}{2}\left(e^{n\delta/\sqrt{p}} + e^{-n\delta/\sqrt{p}}\right).$$

Because $(n+p)\delta^2 = o(n\delta)$ and $n\delta \to \infty$, the right hand side tends to $\infty$. The proof is thus complete. □

*Proof of Lemma 3.3.* First, one gets

$$K(p, n, \delta) = e^{-(n+p)\delta^2/2} \sum_{k=0}^{\infty} \prod_{j=0}^{k-1} \left[\frac{1+2j/(n+p)}{1+2j/p}\right] \times \frac{1}{k!}\left[\frac{(n+p)^2\delta^2}{2p}\right]^k$$

$$\leq e^{-(n+p)\delta^2/2} \sum_{k=0}^{\infty} \frac{1}{k!}\left[\frac{(n+p)^2\delta^2}{2p}\right]^k$$

$$= e^{(n+p)^2\delta^2/(2p) - (n+p)\delta^2/2} = \exp\left[\frac{n(n+p)\delta^2}{2p}\right].$$

Thus, by dominated convergence, $K(p, n, \delta) \to e^a$ as $n(n+p)\delta^2/(2p) \to a$.

Now let $a > 0$. Regard $f(n) = n(n+p)\delta^2/(2p)$ as a quadratic function of $n$. Then in order to get $f(n) \to a$,

$$n \sim \frac{-\delta^2 p + \sqrt{\delta^4 p^2 + 8\delta^2 pa}}{2\delta^2} = \frac{4pa}{\delta^2 p + \sqrt{\delta^4 p^2 + 8\delta^2 pa}}$$

$$\sim \begin{cases} (1/\delta)\sqrt{2pa} & \text{if } \delta^2 p \to 0, \\ 2a/\delta^2 & \text{if } \delta^2 p \to \infty, \\ \dfrac{4a/\delta^2}{1 + \sqrt{1 + 8a/L}} & \text{if } \delta^2 p \to L > 0. \end{cases}$$

The proof is thus complete. □

In order to prove Lemma 3.4, we need the following result.

**Lemma A2.1.** *Given $0 < \epsilon < 1$, there is $\lambda(\epsilon) > 0$, such that*

$$\sum_{|k-A| \geq \epsilon A} \frac{A^k}{k!} \leq e^{A(1-\lambda(\epsilon))}, \qquad as \quad A \to \infty.$$

*Proof.* Let $Y$ be a Poisson random variable with mean $A$. Then

$$e^{-A} \sum_{|k-A| \geq \epsilon A} \frac{A^k}{k!} = P(|Y - A| \geq \epsilon A).$$

By LDP (5), $I := -(1/A) \ln P(|Y - A| \geq \epsilon A) > 0$. Then given $\lambda(\epsilon) \in (0, I)$, $P(|Y - A| \geq \epsilon A) \leq e^{-\lambda(\epsilon)A}$ for all $A \gg 0$, implying the stated bound. □

*Proof of Lemma 3.4.* Fix $\delta > 0$ and $n$. Then

$$K(p, n, \delta) = e^{-A} \sum_{k=0}^{\infty} \prod_{j=0}^{k-1} \left( \frac{n + p + 2j}{p + 2j} \right) \times \frac{A^k}{k!}, \quad \text{with} \quad A = \frac{(n+p)\delta^2}{2}.$$

Let $0 < \epsilon < 1$. For each $k$, $\prod_{j=0}^{k-1}[(n+p+2j)/(p+2j)] \leq (1+n/p)^k$. Then

$$e^{-A} \sum_{|k-A| \geq \epsilon A} \prod_{j=0}^{k-1} \left( \frac{n + p + 2j}{p + 2j} \right) \times \frac{A^k}{k!} \leq e^{-A} \sum_{|k-A| \geq \epsilon A} \frac{[(1+n/p)A]^k}{k!}$$

Denote $B = (1 + n/p)A$. Then given any $0 < \delta < \epsilon$, for all $p \gg 1$, $|k - A| \geq \epsilon A$ implies $|k - B| \geq \delta B$. By Lemma A2.1, as $p \to \infty$,

$$e^{-A} \sum_{|k-A| \geq \epsilon A} \frac{[(1+n/p)A]^k}{k!} \leq e^{-A} \sum_{|k-B| \geq \epsilon B} \frac{B^k}{k!} \leq e^{B - \lambda(\delta)B - A} = o(1),$$

where $\lambda(\delta) > 0$ is a constant. It follows that

$$K(p, n, \delta) = e^{-A} \sum_{|k-A| \leq \epsilon A} \prod_{j=0}^{k-1} \left( 1 + \frac{n}{p + 2j} \right) \times \frac{A^k}{k!} + o(1).$$

By $\ln(1+x) = x + O(x^2)$ as $x \to 0$, it is seen that

$$K(p, n, \delta) = e^{-A} \sum_{|k-A| \leq \epsilon A} (1 + r_k) \exp\left( \frac{n}{p} \sum_{j=0}^{k-1} \frac{1}{1 + 2j/p} \right) \times \frac{A^k}{k!} + o(1),$$

where $\sup_{|k-A| \leq \epsilon A} |r_k| \to 0$ as $p \to \infty$. It is not hard to see that for all $p \gg 1$ and $k$ with $|k - A| \leq \epsilon A$, $|k/p - \delta^2/2| \leq \epsilon \delta^2$. As a result,

$$(1 + r_k) \exp\left( \frac{n}{p} \sum_{j=0}^{k-1} \frac{1}{1 + 2j/p} \right) = [1 + r'_k(\epsilon)] \exp\left( n \int_0^{\delta^2/2} \frac{dx}{1 + 2x} \right)$$

$$= [1 + r'_k(\epsilon)](1 + \delta^2)^{n/2}.$$

where $\sup_{|k-A|\leq \epsilon A} |r'_k(\epsilon)| \to 0$ as $p \to \infty$ followed by $\epsilon \to 0$. Combining the above approximations and applying Lemma A2.1 again,

$$K(p,n,\delta) = [1+R(\epsilon)](1+\delta^2)^{n/2} e^{-A} \sum_{|k-A|\leq \epsilon A} \frac{A^k}{k!} + o(1)$$
$$= [1+R(\epsilon)](1+\delta^2)^{n/2} + o(1),$$

where $R(\epsilon) \to 0$ as $p \to \infty$ followed by $\epsilon \to 0$. Let $p \to \infty$. Since $\epsilon$ is arbitrary, then $K(p,n,\delta) \to (1+\delta^2)^{n/2}$. □

## A3. General $t$ tests

### A3.1. Proof of the main result

This section is devoted to the proof of Proposition 4.1. Write

$$\begin{aligned}\Lambda^*(u) = \sup_t[ut - \Lambda(t)], &\quad \Psi^*(u) = \sup_t[ut - \Psi(t)], \\ \eta_\Lambda(u) = (\Lambda')^{-1}(u), &\quad \eta_\Psi(u) = (\Psi')^{-1}(u),\end{aligned} \quad (A3.1)$$

whenever the functions are well defined. The lemma below collects some useful properties of $\Lambda$. The proof is standard and hence omitted for brevity.

**Lemma A3.1.** *Suppose condition a) in Theorem 4.1 is fulfilled. Then the following statements on $\Lambda$ are true.*

*1) $\Lambda$ is smooth on $\mathcal{D}_\Lambda^o$, strictly decreasing on $(-\infty, 0) \cap \mathcal{D}_\Lambda$, strictly increasing on $(0, \infty) \cap \mathcal{D}_\Lambda$.*

*2) $\Lambda'$ is strictly increasing on $\mathcal{D}_\Lambda^o$, and so $\eta_\Lambda = (\Lambda')^{-1}$ for well defined on $I_\Lambda = (\inf \Lambda', \sup \Lambda')$, where the extrema are obtained over $\mathcal{D}_\Lambda^o$. Moreover, $\Lambda'(0) = 0$, $(\Lambda')^{-1}(0) = 0$, and $t\Lambda'(t) \to \infty$ as $t \uparrow \sup \mathcal{D}_\Lambda$.*

*3) $\Lambda^*$ is smooth and strictly convex on $I_\Lambda$, and*

$$(\Lambda^*)'(u) = \eta_\Lambda(u) = \arg\sup_t[ut - \Lambda(t)], \qquad u \in I_\Lambda.$$

*On the other hand, $\Lambda^*(u) = \infty$ on $(-\infty, \inf \Lambda') \cup (\sup \Lambda', \infty)$.*

The next lemma is key to the proof of Proposition 4.1. Basically, it says that the analysis on the ratio of the extreme tail probabilities can be localized around a specific value determined by $\Lambda$ and the index $\lambda$ in (4.4). As a result, the limit (4.7) can be obtained by the uniform exact large deviations principle (LDP) in (2), which is a refined version of the exact LDP (5).

**Lemma A3.2.** *Let $m, n \to \infty$, such that $n/N \to \rho \in (0,1)$, where $N = m + n$. Let $\nu_0 = \Lambda'(t_0)$, where $t_0 > 0$ the unique positive solution to (4.6). Under conditions a) and b) of Theorem 4.1, given $D > 0$ and $\delta > 0$, there are $z_0 > 0$ and $\eta > 0$, such that for $z \geq z_0$,*

$$\varliminf_{N \to \infty} \frac{1}{N} \inf_{|s| \leq D/N} \ln P\left(\bar X_n + s \geq zS_m,\ |zS_m - \nu_0| \leq \delta\right) \geq -J_z(\nu_0) \quad \text{(A3.2)}$$

*and*

$$\sup_{|s| \leq D/N} \left| \frac{P\left(\bar X_n + s \geq zS_m\right)}{P\left(\bar X_n + s \geq zS_m,\ |zS_m - \nu_0| \leq \delta\right)} - 1 \right| = O(e^{-\eta N}), \quad \text{(A3.3)}$$

*where $J_z(\nu_0) = (1-\rho)\Lambda^*(\nu_0) - \rho\Psi^*(\nu_0^2/z^2) < \infty$.*

Assume Lemma A3.2 is true for now. The main result is shown next.

*Proof of Proposition 4.1.* Recall that $d_N \to 0$ and $N \to \infty$, such that $d_N N \to T$. First, we show that, given $\epsilon > 0$, there is $z_0 > 0$, such that

$$\varlimsup_{N \to \infty} \left| \frac{P\left(\bar X_n + d_N \geq zS_m\right)}{P\left(\bar X_n \geq zS_m\right)} - e^{(1-\rho)Tt_0} \right| \leq \epsilon, \qquad \text{all } z \geq z_0. \quad \text{(A3.4)}$$

Let $\delta \in (0,1)$ such that $\eta_\Lambda(u)$ is well defined on $[\nu_0 - \delta, \nu_0 + \delta]$ and

$$\sup_{|u - \nu_0| \leq \delta} |\eta_\Lambda(u) - \eta_\Lambda(\nu_0)| \leq \frac{\ln(1+\epsilon)}{(1-\rho)T}.$$

Let $z_0 > 0$ and $\eta > 0$ such that (A3.3) holds. Fix $z \geq z_0$. Denote $a = a(z) = (\nu_0 - \delta)/z$ and $b = b(z) = (\nu_0 + \delta)/z$. Because of (A3.3), in order to show (A3.4), it suffices to establish

$$\varlimsup_{N \to \infty} \left| \frac{P\left(\bar X_n + d_N \geq zS_m,\ a \leq S_m \leq b\right)}{P\left(\bar X_n \geq zS_m,\ a \leq S_m \leq b\right)} - e^{(1-\rho)Tt_0} \right| \leq \epsilon. \quad \text{(A3.5)}$$

Let $G_m(x)$ be the distribution function of $S_m$. Then

$$P\left(\bar X_n + d_N \geq zS_m,\ a \leq S_m \leq b\right) = \int_a^b P(\bar X_n \geq zx - d_N)\, G_m(dx),$$

$$P\left(\bar X_n \geq zS_m,\ a \leq S_m \leq b\right) = \int_a^b P(\bar X_n \geq zx)\, dG_m(x).$$

From these equations, it is not hard to see that (A3.5) follows if we can show

$$\varlimsup_{N \to \infty} \sup_{x \in [a,b]} \left| \frac{P(\bar X_n \geq zx - d_N)}{P(\bar X_n \geq zx) e^{(1-\rho)Tt_0}} - 1 \right| \leq \epsilon. \quad \text{(A3.6)}$$

To establish (A3.6), observe that for $N > 1$ large enough and $x \in [a, b]$, $zx - d_N \in [a/2, \nu_0 + \delta]$. Therefore, $\tau_N(x) := \eta_\Lambda(zx - d_N)$ is not only well defined but also continuous and strictly positive on $[a,b]$. By Theorem 3.3 of (2), as $N \to \infty$, the following approximation holds,

$$\sup_{x \in [a,b]} \left| e^{n\Lambda^*(zx-d_N)} \tau_N(x) \sqrt{2\pi n \Lambda''(\tau_N(x))} P(\bar{X}_n \geq zx - d_N) - 1 \right| = o(1),$$

which is a uniform version of the exact LDP due to Bahadur and Rao (5, Theorem 3.7.4).

Because $\tau_N(x) \to \eta_\Lambda(zx)$ uniformly on $[a,b]$ and the latter is strictly positive and continuous on $[a,b]$, the above inequality yields

$$\sup_{x \in [a,b]} \left| e^{n\Lambda^*(zx-d_N)} \eta_\Lambda(zx) \sqrt{2\pi n \Lambda''(\eta_\Lambda(zx))} P(\bar{X}_n \geq zx - d_N) - 1 \right| = o(1).$$

Likewise,

$$\sup_{x \in [a,b]} \left| e^{n\Lambda^*(zx)} \eta_\Lambda(zx) \sqrt{2\pi n \Lambda''(\eta_\Lambda(zx))} P(\bar{X}_n \geq zx) - 1 \right| = o(1).$$

By the above approximations to $P(\bar{X}_n \geq zx - d_N)$ and $P(\bar{X}_n \geq zx)$, in order to prove (A3.6), it is enough to show

$$M := \overline{\lim_{N \to \infty}} \sup_{x \in [a,b]} \left| \frac{e^{-n\Lambda^*(zx-d_N)}}{e^{-n\Lambda^*(zx)+(1-\rho)Tt_0}} - 1 \right| \leq \epsilon.$$

By Taylor expansion and Lemma A3.1,

$$\Lambda^*(zx - d_N) = \Lambda^*(zx) - d_N \eta_\Lambda(zx - \xi d_N), \qquad x \in [a,b],$$

where $\xi = \xi(x) \in (0,1)$. Therefore,

$$\frac{e^{-n\Lambda^*(zx-d_N)}}{e^{-n\Lambda^*(zx)+(1-\rho)Tt_0}} = \frac{e^{-n(\Lambda^*(zx-d_N)-\Lambda^*(zx))}}{e^{(1-\rho)Tt_0}} = \frac{e^{-nd_N \eta_\Lambda(zx-\xi d_N)}}{e^{(1-\rho)Tt_0}}.$$

Since $n d_N \to (1-\rho)T$ and $\eta_\Lambda(zx - \xi d_N) \to \eta_\Lambda(zx)$ uniformly for $x \in [a,b]$,

$$M = \sup_{x \in [a,b]} \left| e^{(1-\rho)(\eta_\Lambda(zx) - t_0)T} - 1 \right|$$

Because $t_0 = \eta_\Lambda(\nu_0)$ and $zx \in [\nu_0 - \delta, \nu_0 + \delta]$ for $x \in [a,b]$,

$$M \leq \exp\left[ (1-\rho)T \sup_{u \in [-\delta, \delta]} |\eta_\Lambda(\nu_0 + u) - \eta_\Lambda(\nu_0)| \right] - 1 \leq \epsilon.$$

Therefore (A3.5) is proved.

Now that (A3.4) holds for any given $\epsilon > 0$, as long as $z \geq z_0 = z_0(\epsilon)$, with $z_0$ being large enough, by the diagonal argument, we can choose $z_N > 0$ in such as way that $z_N \to \infty$ slowly as $N \to \infty$ and

$$\overline{\lim_{N \to \infty}} \left| \frac{e^{-(1-\rho)Tt_0} P(\bar{X}_n + d_N \geq z_N S_m)}{P(\bar{X}_n \geq z_N S_m)} - 1 \right| = 0.$$

This finishes the proof of the theorem. □

### A3.2. Proof of Lemma A3.2

The proof needs a few preliminary results. The first lemma collects some useful properties of $\Psi$.

**Lemma A3.3.** *Let $\mathcal{D}_\Psi = \{t : \Lambda(t) < \infty\}$. Under condition b) in Theorem 4.1, the following statements on $\Psi$ are true.*

1) $\mathcal{D}_\Psi \supset (-\infty, 0]$. $\Psi$ *is smooth and strictly increasing on $\mathcal{D}_\Psi^o$. Furthermore, $\Psi(t) \to -\infty$ as $t \to -\infty$.*

2) $\Psi'$ *is strictly increasing on $\mathcal{D}_\Psi^o$, and so $\eta_\Psi = (\Psi')^{-1}$ is well defined on $I_\Psi = (0, \sup \Psi')$, where the supremum is obtained over $\mathcal{D}_\Psi^o$. In addition, $\inf \Psi' = 0$ and $\sup \Psi' \geq \Psi'(0-) = \sigma^2$. Furthermore,*

$$\lim_{u \to 0+} u \eta_\Psi(u) = -\frac{\lambda + 1}{2}, \qquad (A3.7)$$

*where $\lambda$ is given in (4.4).*

3) $\Psi^*$ *is smooth and strictly convex on $I_\Psi$ and*

$$(\Psi^*)'(u) = \eta_\Psi(u) = \arg\sup_t[ut - \Psi(t)], \qquad u \in I_\Psi.$$

*Furthermore, $\Psi^*$ is strictly decreasing on $(0, \sigma^2)$ with $\Psi^*(u) \to \infty$ as $u \downarrow 0$ and $\Psi^*(u) \to 0$ as $u \uparrow \sigma^2$, and is nondecreasing for $u \geq \sigma^2$.*

*Proof.* We only show $\Psi(t) \to -\infty$ as $t \to -\infty$ and (A3.7), which are properties specifically due to condition b) in Theorem 4.1. The proof of the rest of Lemma A3.3 is standard.

To get $\Psi(t) \to -\infty$ as $t \to \infty$, it suffices to show $\int_0^\infty e^{-tu^2/2} g(u)\, du \to 0$ as $t \to \infty$. For later use, it will be shown that, given $s \geq 0$,

$$\int_0^\infty x^s e^{-tx^2/2} g(x)\, dx \to 0, \qquad \text{as } t \to \infty. \qquad (A3.8)$$

The proof is based on several truncations of the integral. Given $0 < \eta < 1$, there is $0 < \epsilon < 1$, such that

$$1 - \eta \leq \frac{g(x)}{x^\lambda \zeta(1/x)} \leq 1 + \eta, \qquad x \in (0, \epsilon).$$

Since $M_\epsilon = \sup_{|x| \geq \epsilon} g(x) < \infty$, given $s \geq 0$, as $t \to \infty$,

$$\int_\epsilon^\infty x^s e^{-tx^2/2} g(x)\, dx \leq e^{-t\epsilon^2/4} M_\epsilon \int_\epsilon^\infty x^s e^{-tx^2/4}\, dx = o(e^{-t\epsilon^2/4}).$$

On the other hand,

$$\int_0^\epsilon x^s e^{-tx^2/2} g(x)\,dx \geq (1-\eta) \int_0^\epsilon x^{s+\lambda} e^{-tx^2/2} \zeta(1/x)\,dx$$

$$\geq (1-\eta)\zeta(1/\epsilon) \int_0^\epsilon x^{s+\lambda} e^{-tx^2/2}\,dx.$$

The right hand side is of the same order as $\int_0^\infty x^{s+\lambda} e^{-tx^2/2}\,dx$, which in turn is of the same order as $t^{-(\lambda+s+1)/2}$. As a result,

$$\int_0^\infty x^s e^{-tx^2/2} g(x)\,dx = (1+o(1)) \int_0^\epsilon x^s e^{-tx^2/2} g(x)\,dx, \qquad \text{as } t \to \infty.$$

Since $g(x)/[x^\lambda \zeta(1/x)] - 1 \in [-\eta, \eta]$ for $x \in (0, \epsilon)$ and $\eta$ is arbitrary, it is seen that in order to prove (A3.8), it suffices to show

$$\int_0^\epsilon x^{s+\lambda} e^{-tx^2/2} \zeta(1/x)\,dx \to 0, \qquad \text{as } t \to \infty. \tag{A3.9}$$

Let $a = \epsilon^2/2$ and $\phi(x) = \zeta(\sqrt{x/2})$. By variable substitution $x = \sqrt{2u/t}$,

$$\int_0^\epsilon x^{s+\lambda} e^{-tx^2/2} \zeta(1/x)\,dx = 2^p t^{-(p+1)} \int_0^{ta} u^p e^{-u} \phi(t/u)\,du, \tag{A3.10}$$

where $p = (s+\lambda-1)/2 > -1$. Therefore, (A3.9) will follow if

$$t^{-(p+1)} \int_0^{ta} u^p e^{-u} \phi(t/u)\,du \to 0, \qquad \text{as } t \to \infty, \tag{A3.11}$$

Note that $\phi(x)$ is increasing and since $u^p e^{-u}$ is integrable, there is $M > 1$, such that $\int_M^\infty u^p e^{-u}\,du \leq \eta \int_0^M u^p e^{-u}\,du$. Then

$$\int_M^{ta} u^p e^{-u} \phi(t/u)\,du \leq \phi(t/M) \int_M^\infty u^p e^{-u}\,du$$

$$\leq \eta \int_0^M u^p e^{-u} \phi(t/u)\,du. \tag{A3.12}$$

Fix $\delta \in (0,1)$ such that $\delta^{p+1} < \eta(1 - \eta^{p+1}\eta)$. Then

$$\int_0^\delta u^p e^{-u} \phi(t/u)\,du = \sum_{k=1}^\infty \int_{\delta^{k+1}}^{\delta^k} u^p e^{-u} \phi(t/u)\,du$$

$$= \sum_{k=1}^\infty \delta^{(p+1)k} \int_\delta^1 u^p \phi\left(\frac{t}{\delta^k u}\right) du.$$

Note that $\phi(t)$ is slowly varying at $\infty$. For $t$ large enough, $\phi(t/u) \leq \eta\phi(t/(\delta u))$ for $u \in [\eta, 1]$. By induction, $\phi(t/(\delta^k u)) \leq \eta^{k-1} \phi(t/u)$, $k \geq 1$. Consequently, by the selection of $\delta$ and the above infinite sum,

$$\int_0^\delta u^p e^{-u} \phi(t/u)\, du \leq \sum_{k=1}^\infty \delta^{(p+1)k} \eta^{k-1} \int_\delta^1 u^p \phi(t/u)\, du$$
$$= \frac{\delta^{p+1}}{1-\delta^{p+1}\eta} \int_\delta^1 u^p \phi(t/u)\, du \leq \eta \int_\delta^M u^p e^{-u} \phi(t/u)\, du. \tag{A3.13}$$

Now given $0 < \delta < M < \infty$, as $\phi$ is increasing and slowly varying at $\infty$,

$$\inf_{\delta \leq u \leq M} \frac{\phi(t/u)}{\phi(t)} = \frac{\phi(t/M)}{\phi(t)} \to 1, \quad \sup_{\delta \leq u \leq M} \frac{\phi(t/u)}{\phi(t)} = \frac{\phi(t/\delta)}{\phi(t)} \to 1.$$

Therefore,

$$\int_\delta^M u^p e^{-u} \phi(t/u)\, du = (1 + o(1))\phi(t) \int_\delta^M u^p e^{-u}\, du, \quad \text{as } t \to \infty. \tag{A3.14}$$

Combine (A3.12) – (A3.14) and note $\delta$ and $M$ are arbitrary. Then

$$\int_0^{ta} u^p e^{-u} \phi(t/u)\, du = (1 + o(1))\phi(t) \int_0^\infty u^p e^{-u}\, du$$
$$= (1 + o(1))\phi(t)\Gamma(p+1), \quad \text{as } t \to \infty. \tag{A3.15}$$

Note $\phi(t) = o(t^{p+1})$ as $t \to \infty$. Therefore, (A3.11) is proved.

Next we prove (A3.7). For $u > 0$ small enough, $\eta_\Psi(u)$ is well defined. Let $t = -\eta_\Psi(u)$. Then $u = \Psi'(-t)$ and $t \to \infty$ as $u \downarrow 0$. Therefore, it suffice to demonstrate $t\Psi'(-t) \to (\lambda+1)/2$, as $t \to \infty$. It is easy to see

$$\Psi'(-t) = \frac{1}{2} \int_0^\infty x^2 e^{-tx^2/2} g(x)\, dx \Big/ \int_0^\infty e^{-tx^2/2} g(x)\, dx, \quad \text{for } t > 0.$$

Following the argument leading to (A3.9), it suffices to show that, given $\lambda \geq 0$,

$$\int_0^\epsilon x^\lambda e^{-tx^2/2} \zeta(1/x)\, dx = \frac{(1+o(1))t}{\lambda+1} \int_0^\epsilon x^{2+\lambda} e^{-tx^2/2} \zeta(1/x)\, dx$$

as $t \to \infty$. Denoting $p = (\lambda+1)/2$, by (A3.10), the above limit will follow if

$$\int_0^{ta} u^{p-1} e^{-u} \phi(t/u)\, du = \frac{1+o(1)}{p} \int_0^{ta} u^p e^{-u} \phi(t/u)\, du, \quad t \to \infty.$$

However, this is implied by (A3.15) and $\Gamma(p+1) = p\Gamma(p)$. $\square$

**Lemma A3.4.** *Given $\rho \in (0,1)$, let $\nu_0 = \Lambda'(t_0)$, where $t_0 > 0$ is the positive solution to (4.6). Then under conditions a) and b) of Theorem 4.1, for any $\delta \in (0, \nu_0)$, there are $z_0 > 0$ and $a > 0$, such that for $z \geq z_0$,*

$$\inf_{|u-\nu_0| \geq \delta} \left\{ (1-\rho)\Lambda^*(u) + \rho\Psi^*(u^2/z^2) \right\} \geq (1-\rho)\Lambda^*(\nu_0) + \rho\Psi^*(\nu_0^2/z^2) + a.$$

*Proof.* The infimum on the left hand side increases as $\delta$ decreases. Since $\nu_0 < \sup \Lambda'$, without loss of generality, let $\delta < \sup \Lambda' - \nu_0$. For $z > 0$, write

$$H_z(u) = (1-\rho)\Lambda^*(u) + \rho\Psi^*(u^2/z^2)$$

Then by Lemma A3.1, for $u \in (0, \sigma^2 z^2) \cap (0, \sup \Lambda')$,

$$H_z'(u) = (1-\rho)\eta_\Lambda(u) + \frac{2\rho u}{z^2}\eta_\Psi(u^2/z^2) \tag{A3.16}$$

For any $\eta \in (0, \nu_0 - \delta)$ and $M \in (\nu_0 + \delta, \sup \Lambda')$, by (A3.7), as $z \to \infty$,

$$uH_z'(u) \to h(u) := (1-\rho)u\eta_\Lambda(u) - \rho(\lambda+1), \qquad \text{uniformly on } [\eta, M].$$

Since $h$ is strictly increasing on $[0, \infty)$, $\nu_0$ is the only positive solution to $h(u) = 0$. Therefore, there is $a_0 > 0$, such that

$$\inf_{u - \nu_0 \geq \delta/2} h(u) \geq a_0, \qquad \sup_{u - \nu_0 \leq -\delta/2} h(u) \leq -a_0.$$

Let $a = (a_0/2)\min\left\{\ln\frac{\nu_0+\delta}{\nu_0+\delta/2}, \ln\frac{\nu_0-\delta/2}{\nu_0-\delta}\right\}$. As $z \to \infty$, $H_z'(u) \to h(u)/u$ uniformly on $[\eta, M]$. Since $h(u) \geq 0$ for $u \in [\nu_0, M]$, and $h(u)/u \geq a_0/M$ for $u \in [\nu_0 + \delta, M]$, it can be seen that for all $z > 0$ large enough and $u \in [\nu_0 + \delta, M]$,

$$H_z(u) - H_z(\nu_0) = \int_{\nu_0}^u H_z'(s)\,ds \geq \frac{1}{2}\int_{\nu_0+\delta/2}^u \frac{h(s)}{s}\,ds \geq \frac{a_0}{2}\int_{\nu_0+\delta/2}^u \frac{ds}{s} \geq a.$$

Likewise, for all $z > 0$ large enough and $u \in [\eta, \nu_0 - \delta]$,

$$H_z(u) - H_z(\nu_0) = \int_u^{\nu_0}[-H_z'(s)]\,ds \geq \frac{a_0}{2}\int_u^{\nu_0-\delta/2}\frac{ds}{s} \geq a.$$

To finish the proof, it suffices to show that there are $M \in (\nu_0, \sup \Lambda')$ and $\eta \in (0, \nu_0)$, such that for all $z > 0$ large enough, $H_z(u)$ is strictly increasing on $(M, \infty)$ and strictly decreasing on $(0, \eta)$.

First, given $z > 0$ large enough, by Lemma A3.3, $H_z(u)$ is increasing for $u \geq z\sigma^2$ and equal to $\infty$ for $u > \sup \Lambda'$. As a result, it is only necessary to consider $u < M' := \min(\sup \Lambda', z\sigma^2)$. Note that if $\sup \Lambda' < \infty$, then for all $z > 0$ large enough, $M' = \sup \Lambda'$; whereas if $\sup \Lambda' = \infty$, $M' \equiv z\sigma^2$.

Let $\varphi(u) = (u^2/z^2)\eta_\Psi(u^2/z^2)$. For $u \in (\nu_0, M')$, $0 \geq \varphi(u) \geq C := \inf_{0 < u < \sigma^2}[u\eta_\Psi(u)] > -\infty$. By Lemma A3.1, there is $\nu_0 < M < \sup \Lambda'$ such that $(1-\rho)M\eta_\Lambda(M) > -2\rho C + 1$. Then by (A3.16) and the fact that $u\eta_\Lambda(u)$ is strictly increasing for $u \in (0, \sup \Lambda')$, $H_z'(u) > 1/u > 0$ for $u \in (M, M')$. Then $H_z$ is strictly increasing on $(M, M')$.

Second, as $u \downarrow 0$, $u\eta_\Lambda(u) \to 0$ and $u^2\eta_\Psi(u^2) \to -(\lambda+1)/2 < 0$. Therefore, by (A3.16), there is $\delta \in (0, \nu_0)$, such that for all $z > 0$ large enough and $u \in (0, \delta)$, $uH_z'(u) \leq -\rho(\lambda+1)/4$. Then $H_z'(u) < 0$ for $u \leq \delta$ and hence $H_z(u)$ is strictly decreasing. This finishes the proof. $\square$

*Proof of Lemma A3.2.* Since the left hand side of (A3.2) is increasing in $\delta$, without loss of generality, assume $\delta \in (0, \nu_0)$. Let $z_0 > \sigma^2/(\nu_0 + \delta)$. Given $z \geq z_0$ and $\epsilon \in (0, \delta)$, for $N > D/\epsilon$ and $s \in [-D/N, D/N] \subset (-\epsilon, \epsilon)$,

$$P\left(\bar{X}_n + s \geq zS_m, \ |zS_m - \nu_0| \leq \delta\right)$$
$$\geq P\left(\bar{X}_n + s \geq zS_m, \ |zS_m - \nu_0| \leq \epsilon\right)$$
$$\geq P\left(\bar{X}_n \geq \nu_0 + 2\epsilon, \ |zS_m - \nu_0| \leq \epsilon\right)$$
$$= P\left(\bar{X}_n \geq \nu_0 + 2\epsilon\right) P\left(\nu_0 - \epsilon \leq zS_m \leq \nu_0 + \epsilon\right).$$

Observe that for $0 \leq a < b$, $a \leq zS_m \leq b$ is equivalent to $ma^2/z^2 \leq \sum_{k=1}^m (Y_{2k-1} - Y_{2k})^2/2 \leq mb^2/z^2$. Also, $\Lambda^*(t)$ is increasing on $(0, \infty)$, $\Psi^*(t)$ is decreasing on $(0, \sigma^2)$, and $(\nu_0 + \epsilon)^2/z^2 < \sigma^2$. Therefore by LDP,

$$\varliminf_{N \to \infty} \frac{1}{N} \inf_{|s| \leq D/N} P\left(\bar{X}_n + s \geq zS_m, \ |zS_m - \nu_0| \leq \delta\right)$$
$$\geq \varliminf_{N \to \infty} \frac{1}{N} \ln P\left(\bar{X}_n \geq \nu_0 + 2\epsilon\right) + \varliminf_{N \to \infty} \frac{1}{N} \ln P\left(|zS_m - \nu_0| \leq \epsilon\right)$$
$$= (1-\rho)\Lambda^*(\nu_0 + 2\epsilon) + \rho\Psi^*((\nu_0 + \epsilon)^2/z^2).$$

Because $\epsilon$ is arbitrary and $\Lambda^*$ and $\Psi^*$ are continuous, (A3.2) is proved.

Consider (A3.3) now. By Lemma A3.4, there is $\eta > 0$, such that for all $z \geq z_0$ and $u \in [0, \nu_0 - \delta/2] \cup [\nu_0 + \delta/2, \infty)$,

$$(1-\rho)\Lambda^*(u) + \rho\Psi^*(u^2/z^2) \geq J_z(\nu_0) + 2\eta. \qquad (A3.17)$$

Let

$$R_- = \sup_{|s| \leq D/N} P\left(\bar{X}_n + s \geq zS_m, \ zS_m \leq \nu_0 - \delta\right),$$
$$R_+ = \sup_{|s| \leq D/N} P\left(\bar{X}_n + s \geq zS_m, \ zS_m \geq \nu_0 + \delta\right)$$

Since the left hand side of (A3.3) is no greater than

$$\frac{R_- + R_+}{\inf_{|s| \leq D/N} P\left(\bar{X}_n + s \geq zS_m, \ |zS_m - \nu_0| \leq \delta\right)},$$

by (A3.2), in order to establish (A3.3), it suffices to show that for $z \geq z_0$,

$$\varliminf_{N \to \infty} \frac{1}{N} \ln \frac{1}{R_\pm} \geq J_z(\nu_0) + \eta.$$

For any $0 < u \leq \nu_0 - \delta$, by (A3.17), there is $r = r(u) \in (0, u/3)$, such that

$$(1-\rho)\Lambda^*(u - 2r) + \rho\Psi^*((u+r)^2/z^2) \geq J_z(\nu_0) + \eta.$$

By $\Psi^*(u) \uparrow \infty$ as $u \downarrow 0$, there is $r_0 \in (0, \nu_0 - \delta)$, such that $\rho\Psi^*(r_0^2/z^2) \geq J_z(\nu_0) + \eta$. Because $I = [0, \nu_0 - \delta]$ is compact, one can choose $u_0 = 0$ and $u_1, \ldots, u_p \in I$, such that $I \subset \cup_{i=0}^n [u_i - r_i, u_i + r_i]$, with $r_i = r(u_i)$ for $i \geq 1$.

It can be seen that, for $N > D/\min(\epsilon, r_0, r_1, \ldots, r_p)$, $R_- \leq \sum_{i=0}^{p} A_i$, where $A_0 = P(zS_m \leq r_0)$ and $A_i = P(\bar{X}_n \geq u_i - 2r_i, |zS_m - u_i| \leq r_i)$, $i \geq 1$. For the latter ones, by the choice of $z$ and $r_i$, $u_i - 2r_i > 0$ and $(u_i + r_i)/z < \sigma^2$. Therefore, by the LDP,

$$\lim_{N\to\infty} \frac{1}{N} \ln \frac{1}{A_i} = (1-\rho)\Lambda^*(u_i - 2r_i) + \rho\Psi^*((u_i + r_i)^2/z^2) \geq J_z(\nu_0) + \eta.$$

Similarly, $\lim(1/N)\ln(1/A_0) \geq J_z(\nu) + \eta$. Since there is only a finite number of $A_i$, $\underline{\lim}(1/N)\ln(1/R_-) \geq J_z(\nu_0) + \eta$. Likewise, $\underline{\lim}(1/N)\ln(1/R_+) \geq J_z(\nu_0) + \eta$. The proof is thus complete. □

## A4. Tests involving likelihood

### A4.1. Proof of the main result

This section is devoted to the proof of Proposition 5.1. The proof is based on several lemmas. Henceforth, let $N = m + n$ and $\nu_0 = \Lambda'(t_0)$, where $t_0$ is the positive solution to (4.6). It will be assumed that as $m \to \infty$ and $n \to \infty$, $m/N \to \rho \in (0,1)$, where $\rho$ is fixed.

**Lemma A4.1.** *Let $\delta \in (0, \nu_0/2)$ and $\epsilon > 0$. There are $z_0 > 0$ and $\theta_0 = \theta_0(z)$, such that given $z \geq z_0$, as $m \to \infty$ and $n \to \infty$,*

$$\sup_{\theta \leq \theta_0} \frac{P_\theta(\bar{X}_n \geq zS_m, |zS_m - \nu_0| \geq \delta)}{P_\theta(\bar{X}_n \geq zS_m, |zS_m - \nu_0| \leq \delta)} \to 0, \tag{A4.1}$$

$$\sup_{0 \leq \theta \leq \theta_0} \frac{P_\theta(\bar{X}_n \geq (1+\epsilon)zS_m, |zS_m - \nu_0| \leq \delta)}{P_\theta(\bar{X}_n \geq zS_m, |zS_m - \nu_0| \leq \delta)} \to 0. \tag{A4.2}$$

**Lemma A4.2.** *Let $a_1 > 0$. Under the conditions of Theorem 5.1, for any $\epsilon > 0$, there are $m_0 > 0$ and $\delta > 0$, such that*

$$\sup_{0 < |t| \leq \delta} \left| E_0(e^{a\bar{U}_m} \mid S_m = t) - e^{aK_f} \right| \leq \epsilon e^{aK_f}, \quad m \geq m_0, \ 0 \leq a \leq a_1.$$

*where $E_0$ is expectation under $P_0$, $K_f$ is defined as in Theorem 5.1 and $\bar{U}_m = (1/m)\sum_{i=1}^{m} U_i$ with $U_i = (Y_{2i-1} + Y_{2i})/2$.*

*Proof of Proposition 5.1.* We shall show that for any $b > 0$, there is $z_0 = z_0(b)$, such that for all $z \geq z_0$,

$$\overline{\lim} \left| \frac{P_{\theta_N}(\bar{X}_n \geq zS_m)}{P_0(\bar{X}_n \geq zS_m)} - L \right| \leq b, \tag{A4.3}$$

where $L = \exp\{(1-\rho)T\Lambda'(t_0) + 2\rho T K_f\}$, and the limit is taken as $m \to \infty$, $n \to \infty$ and $\theta_N \to 0$, such that $\theta_N N \to T > 0$ and $m/N \to \rho \in (0,1)$. This together with a diagonal argument then finishes the proof.

Let $\epsilon > 0$ and $\delta \in (0, \nu_0/2)$, such that Lemma A4.2 holds with $a = 2\rho T$. Fix $z_0 > (\nu_0 + \epsilon)/\delta$ as in Lemma A4.1. Then, given $z \geq z_0$, in order to show (A4.3), it is enough to show

$$\overline{\lim} \left| \frac{P_{\theta_N}(\mathcal{E}_{m,n})}{P_0(\mathcal{E}_{m,n})} - L \right| \leq b, \tag{A4.4}$$

where $\mathcal{E}_{m,n} = \{\bar{X}_n \in [zS_m, (1+\epsilon)zS_m], |zS_m - \nu_0| \leq \delta\}$. For $\theta \in [0,1]$,

$$\frac{P_\theta(\mathcal{E}_{m,n})}{P_0(\mathcal{E}_{m,n})} = E_0 \left[ e^{J_n(\theta) + mZ_m(\theta)} \,\middle|\, \mathcal{E}_{m,n} \right]$$

where

$$J_n(\theta) = \sum_{i=1}^n \ln r_\theta(\omega_i), \quad Z_m(\theta) = \frac{1}{m} \sum_{j=1}^m \left[ \ln r_\theta(\omega'_{2j-1}) + \ln r_\theta(\omega'_{2j}) \right].$$

Since $\ln r_\theta(\omega_i) = \ell_\theta(\omega_i) - \ell_0(\omega_i)$ and $\dot{\ell}_\theta(\omega_i) = X_i$, by Taylor's expansion,

$$J_n(\theta) = n\theta \bar{X}_n + \frac{\theta^2}{2} \sum_{i=1}^n \ddot{\ell}_{s\theta}(\omega_i), \qquad \text{for some } s \in (0,1).$$

Let $B = \sup_\theta \|\ddot{\ell}_\theta(\omega)\|_{L^\infty(P_0)}$. By Condition 3, $B < \infty$. Since $\theta_N N \to T$, $n\theta_N \to (1-\rho)T$ and $|J_n(\theta_N) - n\theta_N \bar{X}_n| \leq nB\theta_N^2/2 = O(1/N)$. On $\mathcal{E}_{m,n}$,

$$|\bar{X}_n - \nu_0| \leq |\bar{X}_n - zS_m| + |zS_m - \nu_0| \leq \epsilon zS_m + \delta \leq \epsilon_1 := \epsilon(\nu_0 + \delta) + \delta.$$

It follows that for $m$ and $n$ large enough,

$$|J_n(\theta_N) - (1-\rho)T\nu_0| \leq \epsilon + |n\theta_N - (1-\rho)T| \bar{X}_n + (1-\rho)T |\bar{X}_n - \nu_0|$$
$$\leq \epsilon_2 := \epsilon + \epsilon(\nu_0 + \epsilon_1) + (1-\rho)T\epsilon_1.$$

Denote $Q_N = E_0[e^{mZ_m(\theta_N)} \,|\, \mathcal{E}_{m,n}]$. We obtain

$$e^{(1-\rho)T\nu_0 - \epsilon_2} Q_N \leq \frac{P_{\theta_N}(\mathcal{E}_{m,n})}{P_0(\mathcal{E}_{m,n})} \leq e^{(1-\rho)T\nu_0 + \epsilon_2} Q_N. \tag{A4.5}$$

Let $\mathcal{A}_m = \{|zS_m - \nu_0| \leq \epsilon\}$. Since $\omega_i$ and $\omega'_j$ are independent, then

$$Q_N = E_0[e^{mZ_m(\theta_N)} \,|\, \mathcal{A}_m].$$

Let $\bar{U}_m$ be defined as in Lemma A4.2. By Taylor's expansion,

$$mZ_m(\theta) = \theta \sum_{i=1}^m (Y_{2i-1} + Y_{2i}) + \frac{\theta^2}{2} \sum_{i=1}^m [\ddot{\ell}_{s\theta}(\omega_{2i-1}) + \ddot{\ell}_{s\theta}(\omega_{2i})]$$
$$= 2m\theta \bar{U}_m + \frac{\theta^2}{2} \sum_{i=1}^m [\ddot{\ell}_{s\theta}(\omega_{2i-1}) + \ddot{\ell}_{s\theta}(\omega_{2i})], \qquad \text{some } s \in (0,1).$$

Then for $m$ large enough, $|mZ_m(\theta_N) - 2m\theta_N \bar{U}_m| \le BT^2/m < \epsilon$, yielding $e^{-\epsilon} \le Q_N/E_0(e^{2m\theta_N \bar{U}_m} \mid \mathcal{A}_m) \le e^\epsilon$. On $\mathcal{A}_m$, $S_m \le (\nu_0 + \epsilon)/z_0 \le \delta$, so by Lemma A4.2, $1 - \epsilon \le E_0(e^{2m\theta_N \bar{U}_m} \mid \mathcal{A}_m)/e^{2m\theta_N K_f} \le 1 + \epsilon$. By combining (A4.5), we thus obtain

$$(1-\epsilon)e^{-\epsilon-\epsilon_2+2(m\theta_N-\rho T)} L \le \frac{P_{\theta_N}(\mathcal{E}_{m,n})}{P_0(\mathcal{E}_{m,n})} \le (1+\epsilon)e^{\epsilon+\epsilon_2+2(m\theta_N-\rho T)} L.$$

Because $\epsilon$ and $\epsilon_2$ are arbitrary and $m\theta_N \to \rho T$, (A4.3) is proved. □

### A4.2. Proof of Lemmas

We need the next result to show Lemma A4.1.

**Lemma A4.3.** *Given $a \in (0,1)$ and $\epsilon > 0$, there is $\theta_0 > 0$, such that*

$$\sup_{\theta \le \theta_0} P_\theta(\mathcal{E}) \le P_0(\mathcal{E})^{1-a} e^{k\epsilon}, \quad \inf_{\theta \le \theta_0} P_\theta(\mathcal{E}) \ge P_0(\mathcal{E})^{1/(1-a)} e^{-k\epsilon} \quad (A4.6)$$

*for all $k \ge 1$ large enough and $\mathcal{E} \subset \Omega^k$. Furthermore, let $\mathcal{E}_k \subset \Omega^k$ be events such that $\underline{\lim}(1/k) \ln P_0(\mathcal{E}_k) > \infty$. Then*

$$\lim_{\theta_0 \to 0} \overline{\lim_{k \to \infty}} \frac{1}{k} \sup_{0 \le \theta \le \theta_0} \left| \ln \frac{P_\theta(\mathcal{E}_k)}{P_0(\mathcal{E}_k)} \right| = 0.$$

*Proof.* Given $a \in (0,1)$, let $\theta' = \theta'(a)$ as in Condition 4. Denote $\boldsymbol{\omega} = (\omega_1, \ldots, \omega_k)$. For each $\theta \in [0, \theta']$, $k \ge 1$, and $\mathcal{E} \subset \Omega^k$, by Hölder's inequality,

$$P_\theta(\mathcal{E}) = E_0 \left[ \mathbf{1}\{\boldsymbol{\omega} \in \mathcal{E}\} r_\theta(\omega_1) \ldots r_\theta(\omega_k) \right]$$
$$\le [E_0 \mathbf{1}\{\boldsymbol{\omega} \in \mathcal{E}\}]^{1-a} \left\{ E_0 \left[ r_\theta(\omega_k)^{1/a} \ldots r_\theta(\omega_k)^{1/a} \right] \right\}^a$$
$$= P_0(\mathcal{E})^{1-a} \left\{ E_0 \left[ r_\theta(\omega)^{1/a} \right] \right\}^{ka}.$$

Therefore, given $\theta_0 \in (0, \theta')$,

$$\sup_{\theta \le \theta_0} P_\theta(\mathcal{E}) \le P_0(\mathcal{E})^{1-a} \exp\left\{ ka \ln E_0 \left[ \sup_{\theta \le \theta_0} r_\theta(\omega)^{1/a} \right] \right\}$$

Likewise, letting $q = 1/a - 1$,

$$P_0(\mathcal{E}) \le P_\theta(\mathcal{E})^{1-a} \left\{ E_\theta \left[ r_\theta(\omega)^{-1/a} \right] \right\}^{ka} = P_\theta(\mathcal{E})^{1-a} \left\{ E_0 \left[ r_\theta(\omega)^{-q} \right] \right\}^{ka}.$$

Since $q > 0$, the above bound yields

$$\inf_{\theta \le \theta_0} P_\theta(\mathcal{E}) \ge P_0(\mathcal{E})^{1/(1-a)} \exp\left\{ -\frac{ka}{1-a} \ln E_0 \left[ (\inf_{\theta \le \theta_0} r_\theta(\omega))^{-q} \right] \right\}$$

Under $P_0$, for almost every $\omega \in \Omega$, $p_0(\omega) > 0$ and $p_\theta(\omega)$ is continuous in $\theta$. Let $\theta_0 \to 0$. Then $\sup_{\theta \leq \theta_0} r_\theta(\omega) \to 1$ and $\inf_{\theta \leq \theta_0} r_\theta(\omega) \to 1$. By (5.2) and dominated convergence,

$$\ln E_0 \left[ \sup_{\theta \leq \theta_0} r_\theta(\omega)^{1/a} \right] \to 0, \quad \ln E_0 \left[ (\inf_{\theta \leq \theta_0} r_\theta(\omega))^{-q} \right] \to 0.$$

This implies that for $\theta_0$ small enough, both of the inequalities in (A4.6) hold.

To show the second part of the lemma, for each $n \geq 1$,

$$\frac{1}{k} \ln P_\theta(\mathcal{E}_k) \leq \frac{1-a}{k} \ln P_0(\mathcal{E}_k) + a \ln E_0[r_\theta(\omega)^{1/a}],$$

which yields

$$\frac{1}{k} \sup_{0 \leq \theta \leq \theta_0} \ln \frac{P_\theta(\mathcal{E}_k)}{P_0(\mathcal{E}_k)} \leq a \left\{ -\frac{1}{k} \ln P_0(\mathcal{E}_k) + \ln E_0 \left[ \sup_{\theta \leq \theta_0} r_\theta(\omega)^{1/a} \right] \right\}.$$

Let $k \to \infty$ and take $\overline{\lim}$ on both ends. By the assumption,

$$\overline{\lim_{k \to \infty}} \frac{1}{k} \sup_{0 \leq \theta \leq \theta_0} \ln \frac{P_\theta(\mathcal{E}_k)}{P_0(\mathcal{E}_k)} \leq a \left\{ M + \ln E_0 \left[ \sup_{\theta \leq \theta_0} r_\theta(\omega)^{1/a} \right] \right\},$$

where $M = -\underline{\lim}(1/k) \ln P_0(\mathcal{E}_k) \geq 0$. Likewise, with $q = 1/a - 1 > 0$,

$$\underline{\lim_{k \to \infty}} \frac{1}{k} \inf_{0 \leq \theta \leq \theta_0} \ln \frac{P_\theta(\mathcal{E}_k)}{P_0(\mathcal{E}_0)} \geq -\frac{a}{1-a} \left\{ M + \ln E_0 \left[ (\inf_{\theta \leq \theta_0} r_\theta(\omega))^{-b} \right] \right\}.$$

Thus we get

$$\overline{\lim_{\theta_0 \to 0}} \; \overline{\lim_{k \to \infty}} \sup_{0 \leq \theta \leq \theta_0} \frac{1}{k} \left| \ln \frac{P_\theta(\mathcal{E}_k)}{P_0(\mathcal{E}_k)} \right| \leq \frac{aM}{1-a}.$$

Because $a$ is arbitrary, the lemma is proved. □

It is easy to check that under the assumptions of Theorem 5.1, all the statements in Lemmas A3.1 and A3.3 hold for $\Lambda$ and $\Psi$ defined in (5.6), with $X = \dot{\ell}_0(\omega)$, $Y = \dot{\ell}_0(\omega')$. Therefore, Lemma A3.2 can be applied.

*Proof of Lemma A4.1.* We first show (A4.1). By Lemma A3.2, there is $z_0 > 0$, such that for $z \geq z_0$ and $\theta \in (0, \nu_0/2)$, there is $\eta > 0$, such that

$$\frac{P_0(\mathcal{E}_{n,m} \cap \mathcal{A}_m^c)}{P_0(\mathcal{E}_{n,m} \cap \mathcal{A}_m)} = o(e^{-\eta M}), \tag{A4.7}$$

where $M = n + 2m$, $\mathcal{E}_{n,m} = \{\bar{X}_n \geq zS_m\}$ and $\mathcal{A}_m = \{|zS_m - \nu_0| \leq \delta\}$. Given $\epsilon \in (0,1)$, by Lemma A4.3, there is $\theta_0 > 0$, such that for $\theta \in [0, \theta_0]$ and $m$, $n$

large enough, $P_\theta(\mathcal{E}) \leq P_0(\mathcal{E})^{1-\epsilon} e^{\epsilon M}$ and $P_\theta(\mathcal{E}) \geq P_0(\mathcal{E})^{1/(1-\epsilon)} e^{-\epsilon M}$ for $\mathcal{E} \subset \Omega^M$. Since both $\mathcal{E}_{n,m}$ and $\mathcal{A}_m$ are events in $\Omega^M$, then

$$L_{n,m} := \frac{1}{M} \sup_{0 \leq \theta \leq \theta_0} \ln \frac{P_\theta(\mathcal{E}_{n,m} \cap \mathcal{A}_m^c)}{P_\theta(\mathcal{E}_{n,m} \cap \mathcal{A}_m)}$$

$$\leq \frac{1}{M} \ln \frac{P_0(\mathcal{E}_{n,m} \cap \mathcal{A}_m^c)^{1-\epsilon}}{P_0(\mathcal{E}_{n,m} \cap \mathcal{A}_m)^{1/(1-\epsilon)}} + 2\epsilon$$

$$= \frac{1}{M} \left[ (1-\epsilon) \ln \frac{P_0(\mathcal{E}_{n,m} \cap \mathcal{A}_m^c)}{P_0(\mathcal{E}_{n,m} \cap \mathcal{A}_m)} + \frac{\epsilon(2-\epsilon)}{1-\epsilon} \ln \frac{1}{P_0(\mathcal{E}_{n,m} \cap \mathcal{A}_m)} \right] + 2\epsilon.$$

By equations (A3.2) and (A4.7), there is a finite constant $C > 0$, such that

$$\overline{\lim} L_{n,m} \leq -(1-\epsilon)\eta + \frac{\epsilon(2-\epsilon)C}{1-\epsilon} + 2\epsilon,$$

Since $\epsilon$ is arbitrary, $\overline{\lim} L_{n,m} < 0$. This then finishes the proof of (A4.1).

It remains to show (A4.2). First, by the LDP for $\bar{X}_n$ under $P_0$ and an argument similar to the proof of (A4.1), it can be seen that given $r > 0$ and $0 < a < b < \sup_{\mathcal{D}_\Lambda^o} \Lambda')$, there is $\theta_0 > 0$, such that

$$\sup_{0 \leq \theta \leq \theta_0} [P_\theta(\bar{X}_n \geq b) / P_\theta(\bar{X}_n \in [a, a+r])] \to 0, \quad \text{as } n \to \infty. \tag{A4.8}$$

Now let $a \in (0, \epsilon)$ and $\eta \in (0, (\delta/\nu_0) \wedge (a/2))$, so that $(1+\epsilon)(1-\eta) > 1+a$. Denote $\mathcal{E}_m = \{|zS_m - \nu_0| \leq \eta\nu_0\}$ and $\mathcal{A}_m = \{|zS_m - \nu_0| \leq \delta\}$. Then $\mathcal{E}_m \subset \mathcal{A}_m$. By Lemma A3.2, given $z \gg 1$, there is $\theta_0 > 0$, such that

$$\inf_{\theta \leq \theta_0} \frac{P_\theta(\bar{X}_n \geq zS_m, \mathcal{E}_m)}{P_\theta(\bar{X}_n \geq zS_m, \mathcal{A}_m)} \to 1. \tag{A4.9}$$

For $\theta \leq \theta_0$, by the independence of $\bar{X}_n$ and $S_m$ under $P_\theta$,

$$P_\theta(\bar{X}_n \geq (1+\epsilon)zS_m, \mathcal{E}_m) \leq P_\theta(\bar{X}_n \geq (1+\epsilon)(1-\eta)\nu_0, \mathcal{E}_m)$$
$$\leq P_\theta(\bar{X}_n \geq (1+a)\nu_0) P_\theta(\mathcal{E}_m).$$

By $\eta < a/2$, let $\epsilon' \in (0, \epsilon)$, such that $(1+\epsilon')(1+\eta) \leq 1 + a/2$. Let $I = [(1-\eta)\nu_0, (1+\epsilon')(1+\eta)\nu_0]$. It is not hard to find a finite number of nonempty $(b_i, c_i) \subset I$, such that for any $x \in I$, $[x, (1+\epsilon')x]$ contains at least one $(b_i, c_i)$. Then

$$P_\theta(\bar{X}_n \in [zS_m, (1+\epsilon)zS_m], \mathcal{E}_m) \geq P_\theta(\bar{X}_n \in [zS_m, (1+\epsilon')zS_m], \mathcal{E}_m)$$
$$\geq \min_i P_\theta(\bar{X}_n \in [b_i, c_i]) P_\theta(\mathcal{E}_m)$$

Since $c_i \leq (1+a/2)\nu_0$, by the above inequalities and (A4.8),

$$\sup_{\theta \leq \theta_0} \frac{P_\theta(\bar{X}_n \geq (1+\epsilon)zS_m, \mathcal{E}_m)}{P_\theta(\bar{X}_n \geq zS_m, \mathcal{E}_m)} \leq \sup_{\theta \in \theta_0} \frac{P_\theta(\bar{X}_n \geq (1+\epsilon)zS_m, \mathcal{E}_m)}{P_\theta(\bar{X}_n \geq [zS_m, (1+\epsilon)zS_m], \mathcal{E}_m)}$$
$$\leq \max_i \sup_{\theta \in \theta_0} \frac{P_\theta(\bar{X}_n \geq (1+a)\nu_0)}{P_\theta(\bar{X}_n \in [b_i, c_i])} \to 0,$$

yielding

$$\inf_{\theta \leq \theta_0} \frac{P_\theta(\bar X_n \in [zS_m,\, (1+\epsilon)zS_m],\, \bar{\mathcal{E}}_m)}{P_\theta(\bar X_n \geq zS_m,\, \mathcal{E}_m)}$$
$$\geq \inf_{\theta \leq \theta_0} \frac{P_\theta(\bar X_n \in [zS_m,\, (1+\epsilon)zS_m],\, \mathcal{E}_m)}{P_\theta(\bar X_n \geq zS_m,\, \mathcal{E}_m)} \to 1,$$

which, together with (A4.9), implies

$$\inf_{\theta \leq \theta_0} \frac{P_\theta(\bar X_n \in [zS_m,\, (1+\epsilon)zS_m],\, \mathcal{A}_m)}{P_\theta(\bar X_n \geq zS_m,\, \mathcal{A}_m)}$$
$$\geq \inf_{\theta \leq \theta_0} \frac{P_\theta(\bar X_n \in [zS_m,\, (1+\epsilon)zS_m],\, \mathcal{A}_m)}{P_\theta(\bar X_n \geq zS_m,\, \mathcal{E}_m)} \times \inf_{\theta \leq \theta_0} \frac{P_\theta(\bar X_n \geq zS_m,\, \mathcal{E}_m)}{P_\theta(\bar X_n \geq zS_m,\, \mathcal{A}_m)} \to 1$$

and hence (A4.2). □

*Proof of Lemma A4.2.* Let $\epsilon_0 = \epsilon \inf_{0 \leq a \leq a_1} e^{aK_f}$. We have to show that for $\delta > 0$ small enough and $m_0 > 0$ large enough,

$$\sup_{0 < |t| \leq \delta} \left| E_0(e^{a\bar U_m} \mid S_m = t) - e^{aK_f} \right| < \epsilon_0, \quad m \geq m_0,\ 0 \leq a \leq a_1. \quad (A4.10)$$

Let $V_i = (Y_{2i-1} - Y_{2i})/2$. Under $P_0$, $(U_i, V_i)$ are iid with density

$$\frac{P(U \in du, V \in dv)}{du\,dv} = 2f(u+v)f(u-v).$$

Denote $\boldsymbol{\zeta} = (v_1, \ldots, v_m)$ and

$$\phi_v(z) = E_0(e^{zU} \mid V = v).$$

Then

$$E_0(e^{a\bar U_m} \mid S_m = t) = \int \prod_{i=1}^m \phi_{v_i}(a/m) \mathbf{1}\{v_i \neq 0\}\, P_0(d\boldsymbol{\zeta} \mid S_m = t). \quad (A4.11)$$

**Case i: $f$ is bounded** In this case, $g(v) = \int f(u+v)f(u-v)\,du$ is well defined for all $v \in \text{sppt}(V)$, $h_v(u) = f(u+v)f(u-v)/g(v)$ is the conditional density of $U$ given $V = v$ and $\phi_v(z) = \int e^{zu} h_v(u)\,du$. Since $f$ is continuous almost everywhere and bounded, by condition a) of Theorem 5.1, there is $r > 0$ such that $\sup_v \int e^{r|u|} f(u+v)f(u-v)\,du < \infty$, and by dominated convergence, as $v \to 0$, $g(v) \to g(0) = \int f^2 \in (0, \infty)$. It follows that there is $c > 0$, such that $\{\phi_v(z), v \in [-c, c]\}$ is a family of smooth functions of $z \in [-r, r]$ with uniformly continuous and bounded $\phi_v'(z)$ and $\phi_v''(z)$.

Given $\eta > 0$, decrease $c$ if necessary so that

$$\sup_{(v,z) \in I} |\phi_v^{(k)}(z) - \phi_0^{(k)}(z)| \leq \frac{\eta}{3a_1}, \quad k = 1, 2,$$

where $I = [-c, c] \times [-r, r]$. By Taylor's expansion,

$$\phi_v(z) - \phi_0(z) = [\phi'_v(0) - \phi'_0(0)]z + \frac{1}{2}[\phi''_v(\theta z) - \phi''_0(\theta z)]z^2, \quad (v, z) \in I,$$

where $\theta = \theta(v, z) \in (0, 1)$. Then there is $m_0 > 0$, such that for all $m \geq m_0$, $a \in [0, a_1]$ and $v \in I$, one gets $a/m \in [-r, r]$,

$$|\phi_v(a/m) - \phi_0(a/m)| \leq 2\eta/(3m) \leq (\eta/m) \inf_{0 \leq a \leq a_1} \phi_0(a/m)$$

and hence

$$1 - \frac{\eta}{m} \leq \frac{\phi_v(a/m)}{\phi_0(a/m)} \leq 1 + \frac{\eta}{m}, \quad \text{all } a \in [0, a_1]. \tag{A4.12}$$

Given $\delta \in (0, c)$, for $0 < t \leq \delta$, rewrite (A4.11) as

$$E_0(e^{a\bar{U}_m} \mid S_m = t) = \int \prod_{i \in J} \phi_{v_i}(a/m) \prod_{i \notin J} \phi_{v_i}(a/m) \, P(d\boldsymbol{\zeta} \mid S_m = t),$$

where $J = \{i : |v_i| \geq c\}$. Let $s > 0$ and $L > 0$ be as in (5.8). For $m$ large enough, $a/m < s$, $a \in [0, a_1]$. Therefore, by Hölder's inequality, for $i \in J$,

$$\phi_{v_i}(a/m) \leq [\phi_{v_i}(s)]^{a/(sm)} \leq L^{a/(sm)} \exp\left(\frac{aL|v_i|}{sm}\right),$$

$$\phi_{v_i}(a/m) \geq \frac{1}{\phi_{v_i}(-a/m)} \geq L^{-a/(sm)} \exp\left(-\frac{aL|v_i|}{sm}\right).$$

Let $p = |J|/m$. By the above first set of inequalities and Schwartz inequality,

$$\prod_{i \in J} \phi_{v_i}(a/m) \leq L^{ap/s} \exp\left(\frac{aL}{sm} \sum_{i \in J} |v_i|\right) \leq L^{ap/s} \exp\left(\frac{aL}{sm}\sqrt{J \sum v_i^2}\right).$$

Likewise, by the above second set of inequalities and Schwartz inequality,

$$\prod_{i \in J} \phi_{v_i}(a/m) \geq L^{-ap/s} \exp\left(-\frac{aL}{sm}\sqrt{J \sum v_i^2}\right)$$

Since $\{S_m = t\} = \{(1/m) \sum V_i^2 = t^2/4\}$,

$$L^{-ap/s} \exp\left(-\frac{aLt\sqrt{p}}{2s}\right) \leq \prod_{i \in J} \phi_{v_i}(a/m) \leq L^{ap/s} \exp\left(\frac{aLt\sqrt{p}}{2s}\right).$$

Observe that, due to $0 < t \leq \delta$, $S_m = t$ implies $p \leq \delta^2/c^2$. Therefore, as long as $\delta$ is small enough, $ap/s$ is arbitrarily close to 0, and $aLt\sqrt{p}/(2s)$ is uniformly arbitrarily close to 0 for $0 \leq a \leq a_1$ and $0 < t \leq \delta$. Consequently, for each $\boldsymbol{\zeta} \in \{S_m = t\}$, $e^{-\eta} \leq \prod_{i \in J} \phi_{v_i}(a/m) \leq e^{\eta}$.

On the other hand, by (A4.12),

$$e^{-\eta} \leq \left(1 - \frac{\eta}{m}\right)^{m(1-p)} \leq \prod_{i \notin J} \frac{\phi_{v_i}(a/m)}{\phi_0(a/m)} \leq \left(1 + \frac{\eta}{m}\right)^{m(1-p)} \leq e^{\eta}.$$

Thus, $e^{-2\eta}\phi_0(a/m)^{m(1-p)} \leq E_0(e^{a\bar{U}_m} \,|\, S_m = t) \leq e^{2\eta}\phi_0(a/m)^{m(1-p)}$ for all $t \in [-\delta, \delta] \setminus \{0\}$. Since $\eta$ and $p$ are arbitrarily small and $\phi_0(a/m)^m \to e^{aK_f}$ uniformly for $a \in [0, a_1]$ as $m \to \infty$, (A4.10) then follows.

**Case ii: $f$ is symmetric and has a bounded support** In this case $B := \|U\|_{L^\infty(P_0)} < \infty$. By condition c) of Theorem 5.1, the density of $V$ is continuous and bounded on $(\epsilon, \infty)$ for any $\epsilon$. Then $\phi_v(z)$ is well defined for all $z$ and $v \in \text{sppt}(V) \setminus \{0\}$. Since $f$ is symmetric, for $v \in \text{sppt}(V) \setminus \{0\}$,

$$\phi_v'(0) = \int u f(u+v) f(u-v) \, du / g(v) = 0,$$

and so $|\phi_v(a/m) - 1| \leq |\phi_v''(\theta a/m)| \, (a/m)^2$, with $\theta \in (0, 1)$. By

$$\phi_v''(s) = \int u^2 e^{su} f(u+v) f(u-v) \, du / g(v) \leq B^2 e^{|s|B}$$

Then $|\phi_v(a/m) - 1| \leq (a/m)^2 B_1$, where $B_1 = B^2 e^{(a/m)B}$. Then by (A4.11), $[1 - B_1(a/m)^2]^m \leq E_0(e^{a\bar{U}_m} \,|\, S_m = t) \leq [1 + B_1(a/m)^2]^m$, which implies (A4.10). □